\documentclass[preprint]{elsarticle}
\usepackage{amsmath}
\usepackage{amsfonts}
\usepackage{amssymb}
\usepackage{bm,bbm}
\usepackage{algorithm}
\usepackage{algorithmic}
\usepackage[aboveskip=1pt]{subcaption}
\usepackage{caption}
\usepackage{diagbox}
\usepackage{lineno,hyperref}
\usepackage[textheight = 21cm, textwidth = 16cm]{geometry}

\journal{Computers \& Fluids, published 15 August 2022, \url{https://doi.org/10.1016/j.compfluid.2022.105575}}

\begin{document}

\begin{frontmatter}

\title{Comparative study of inner-outer Krylov solvers for linear systems in structured and high-order unstructured CFD problems}

\author[mymainaddress]{Mehdi Jadoui}
\ead{Mehdi.Jadoui@onera.fr}

\author[mymainaddress]{Christophe Blondeau\corref{mycorrespondingauthor}}
\cortext[mycorrespondingauthor]{Corresponding author}
\ead{Christophe.Blondeau@onera.fr}

\author[mymainaddress]{Emeric Martin}
\ead{Emeric.Martin@onera.fr}

\author[mymainaddress]{Florent Renac}
\ead{Florent.Renac@onera.fr}

\author[mysecondaryaddress]{François-Xavier Roux}
\ead{roux@ann.jussieu.fr}

\address[mymainaddress]{ONERA, Université Paris-Saclay, F-92322 Châtillon, France}
\address[mysecondaryaddress]{Laboratory of Jacques-Louis Lions, Sorbonne University, F-75005 Paris, France }

\begin{abstract}
Advanced Krylov subspace methods are investigated for the solution of large sparse linear systems arising from stiff adjoint-based aerodynamic shape optimization problems. A special attention is paid to the flexible inner-outer GMRES strategy combined with most relevant preconditioning and deflation techniques. The choice of this specific class of Krylov solvers for challenging problems is based on its outstanding convergence properties. Typically in our implementation the efficiency of the preconditioner is enhanced with a domain decomposition method with overlapping. However, maintaining the performance of the preconditioner may be challenging since scalability and efficiency of a preconditioning technique are properties often antagonistic to each other. In this paper we demonstrate how flexible inner-outer Krylov methods are able to overcome this critical issue. A numerical study is performed considering either a Finite Volume (FV), or a high-order Discontinuous Galerkin (DG) discretization which affect the arithmetic intensity and memory-bandwith of the algebraic operations. We consider test cases of transonic turbulent flows with RANS modelling over the two-dimensional supercritical OAT15A airfoil and the three-dimensional ONERA M6 wing. Benefits in terms of robustness and convergence compared to standard GMRES solvers are obtained. Strong scalability analysis shows satisfactory results. Based on these representative problems a discussion of the recommended numerical practices is proposed.
\end{abstract}

\begin{keyword}
Krylov solver, GMRES, adjoint problem, deflation, variable preconditioning 
\end{keyword}

\end{frontmatter}


\section{Introduction}

One of the most important environmental challenge that aircraft manufacturers seek to overcome is the limitation of the carbon footprint in the atmosphere. Many technological improvements have been achieved to reduce harmful emissions. One of them concerns the improvement of the efficiency of the airfoil that is typically quantified by the lift to drag ratio. More specifically, one seeks to minimize an aerodynamic function with respect to a set of design variables that controls the shape of the airfoil. Gradient-based methods are the most suitable optimization algorithms to address such problems. The efficient computation of the sensitivities is consequently an essential step to obtain a relevant design. In aerodynamic design the adjoint formulation allows the computation of derivatives of several functions with respect to many design parameters at an affordable cost. However, high-fidelity analyses are now routinely applied even at the preliminary design stage and it becomes increasingly challenging to solve the resulting algebraic systems with regard to the size and the condition number of such problems. At this stage, it is obvious that the accuracy of the gradient computation depends on the accuracy of the solution of the linear system.

In spite of the rise of supercomputers and their high performance computing, solving linear systems is still a substantive issue in several scientific fields like computational fluid dynamics (CFD). Indeed, predicting the best design requires models of compressible turbulent flows in order to capture complex physical phenomena. From an algebraic point of view, this leads to large, sparse and ill-conditioned systems. Matrices are real, non symmetric, not positive definite, with a block-wise structure and a symmetric pattern. A robust way to solve this kind of systems relies on Krylov subspace methods, an iterative procedure which computes an approximation of the solution in a projection subspace called Krylov subspace. For non-symmetric problems, we focus on the Generalized Minimal Residual (GMRES) algorithm~\cite{saad1986}. However, due to the prohibitive memory cost of the GMRES, the restarted version denoted as GMRES($\textit{m}$) is practically used, with $\textit{m}$ the maximum number of basis vectors allowed during a cycle. Once $\textit{m}$ iterations are carried out, the vectors are discarded and a new set of vectors is generated from scratch. The counterpart is the loose of the super-linearity property of the GMRES and therefore a convergence slowdown may be observed compared to the full GMRES without restart. Numerical ingredients as preconditioning techniques and spectral deflation have been therefore introduced to offset the restart effect. Generally speaking, preconditioners attempt to improve the spectral properties of the operator of the system by clustering as much as possible the eigenspectrum of the initial system. Benzi~\cite{benzi2002} surveyed different preconditioning techniques for the iterative solution of large linear systems. For instance, he mentioned the incomplete factorization preconditioner and its block variant denoted by ILU and BILU respectively. The ILU algorithm seems to be attractive thanks to its robustness and as a result is used as a preconditioner in most cases. But some limitations exist about stability and scalability. Stability is in fact affected by many parameters like the condition number, and the structure of the coefficient matrix. For our case, instability issues are avoided since the preconditioning operators have good algebraic properties. The BILU algorithm will be employed throughout this paper. Another class of well-known iterative preconditioners is the Lower-Upper Symmetric Successive Over Relaxation (LU-SSOR) preconditioner combining the LU factorization with an SSOR relaxation~\cite{yoon1988}. However, convergence of iterative techniques may be hampered by the lack of diagonal dominance of the associated operator.

A seminal contribution to efficiency of preconditioned GMRES was the extension to variable preconditioning~\cite{saad1993flexible}. This gives GMRES a great flexibility to use advanced iterative solvers essentially as preconditioners which may themselves in turn consist in Krylov subspace methods. When it comes to use GMRES as a preconditioner, the method is called Flexible inner-outer GMRES or nested GMRES. However, it could be even difficult for the nested GMRES methods to address some extreme ill-conditioned systems. A solution consists in using the spectral deflation strategy. It aims at accelerating the convergence when preconditioned FGMRES fails to converge. This stagnation is essentially caused by the presence of eigenvalues in the neighborhood of zero. The idea is to find a way to remove their harmful effect from the eigenspectrum of the preconditioned system. Morgan~\cite{morgan2002} proposed to deflate the harmonic Ritz vectors from the next Krylov subspace since they are a good approximation of eigenvectors associated with the smallest eigenvalues. The harmonic Ritz vectors are formulated as a linear combination of Krylov vectors whose coefficients are solutions of a generalized eigenvalue problem~\cite{morgan1998harmonic}. It leads to a new algorithm denoted as GMRES-DR($m$,$k$) (GMRES with Deflating Restarting) with $k$ the number of selected harmonic Ritz vectors per cycle. Giraud et al.~\cite{giraud2010} proposed to generalize this algorithm to a flexible formulation, FGMRES-DR($m$,$k$), where the preconditioner changes at each Arnoldi iteration. They found effective to deflate eigenvalues that are located at a distance away from specific clusters. In the field of aerodynamic shape optimization, few research has been conducted regarding the inner-outer GMRES method. In a previous work, Pinel and Montagnac~\cite{pinel2013block} investigated the nested GMRES method as well as its block variant for the solution of a discrete adjoint problem. The latter problem results from the linearization of the discretized Navier-Stokes equations combined with the two-equation turbulence model of Wilcox. We point out that the eddy viscosity is  kept constant during the linearization. They have proved that the block nested GMRES overcomes the scalar version with a factor of 3 in terms of CPU time. In addition to that, the same performance is reached with lower memory footprint for the block nested GMRES. Chen et al.~\cite{chen2019gcro} recently tested both GMRES-DR and GCRO-DR (Generalized Conjugate Residual with inner Orthogonalization and Deflated Restarting)\cite{parks2006recycling} solvers by proposing an extension in terms of dynamic deflated restarting. The results demonstrate the effectiveness of this strategy on adjoint systems in terms of memory cost and CPU time. 
        
This paper aims at investigating the flexible inner-outer GMRES solver capabilities to solve adjoint systems for turbulent flow by fully linearizing the turbulence model. Actually, only a handful of authors have examined this option due to the following reason~\cite{Dwight2006, Nemec2001,Renac2007,Renac2011} (see the review~\cite{Peter2010} and references therein). Linearizing the turbulence model leads generally to very ill-conditioned adjoint systems demanding the development of advanced linear solvers combined with others numerical artifices. However, considering the turbulence model during the linearization is crucial when it comes to solve a high-fidelity optimization problem. So, instead of making physical simplifications, a special focus on nested GMRES solvers associated with spectral deflation and preconditioning techniques has been achieved. In contrast to~\cite{pinel2013block}, a special attention to the block incomplete LU factorization (BILU) is conferred since it is proved to be more robust than the standard block Lower-Upper Symmetric Gauss-Seidel (LU-SGS) relaxation algorithm. LU-SGS is a particular case of LU-SSOR with a unity relaxation factor. The reader should refer to~\citep{peter2007large} for further detail about LU-SGS as an iterative solver for the implicit phase of Backward-Euler Newton steady solvers. Numerical experiments have shown that convergence is strongly affected depending on the choice of the preconditioner. We also point out that we recently implemented the GCRO-DR and a comparative study with GMRES-DR has shown that performances are strictly equivalent. For these reasons, the current work has been focused only on the nested GMRES solver. Considering two different discretization methods, such as FV and DG formulations, contributes to draw a better overview of the benefits of such a solver. Not only it affects conditioning of adjoint systems but also the arithmetic intensity and the memory-bandwidth pressure of basic linear algebra operations.

This paper is organized as follows. In section 2, we briefly present the adjoint problem for a sensitivity analysis. In section 3, we describe the minimal residual norm Krylov subspace methods in conjunction with deflation strategies and flexible preconditioning. Section 4 and 5 are devoted to numerical experiments where the robustness of flexible inner-outer Krylov solvers is demonstrated through representative test-cases for both FV \footnote{FV simulations were performed with the ONERA elsA software~\cite{cambier2013onera}} and DG \footnote{DG simulations were performed with the ONERA Aghora software~\cite{Renac2015}} schemes. Strong scalability capability of such solvers is also assessed. To finis with, some numerical considerations are discussed in section 6.

\section{Discrete adjoint method for sensitivity analysis}

\subsection{Steady-state problem}

Let $\mathcal{D}$ be a bounded domain in $\mathbb{R}^3$. We denote by $\textbf{W} \in \mathbb{R}^{N}$ the vector of conservative variables.
The differential form of the governing equations for a viscous fluid governed by the RANS model is formulated as

\begin{equation}
    \label{eq:RANS_model}
    \frac{d\textbf{W}}{dt} + \nabla . \textbf{F} = \textbf{Q},
\end{equation}
where $\textbf{F}$ embeds the convective and diffusive fluxes and $\textbf{Q}$ is a source term. We want to obtain the field of conservative variables \textbf{W} related to a steady-state equilibrium on domain $\mathcal{D}$. The system of equations is formulated in its discrete residual form as:

\begin{equation}
    \label{eq:Discret_Residual}
    \mathbf{R}(\mathbf{W}, \mathbf{X}) = \mathbf{0},
\end{equation}
where $\mathbf{X} \in \mathbb{R}^N$ is the vector of fluid grid coordinates.

\subsection{Sensitivity analysis via the adjoint approach}

Consider an aerodynamic shape optimization problem of the form $ \mathcal{J}^* = \underset{\alpha}{\min} \,\mathcal{J}(\mathbf{W}( \alpha), \mathbf{X}(\alpha))$ under the constraints $\mathbf{R}(\mathbf{W}(\alpha), \mathbf{X}(\alpha)) = \mathbf{0}$, where $\alpha$ is the vector of design variables controlling the shape of the body. A gradient-based method is required to efficiently minimize $\mathcal{J}$. In aerodynamics, the adjoint technique is the method of choice to compute the sensitivities~\cite{jameson1995optimum,Peter2010}. The adjoint linear system is written as:

\begin{equation}
    \label{eq:Adjoint_System}
    \left[ \frac{ \partial \mathbf{R}}{ \partial \mathbf{W}} \right]^T \lambda  =  - \frac{ \partial \mathcal{J}^{T}}{ \partial \mathbf{W}} 
\end{equation}
where the superscript $^{T}$ stands for the transpose operator. In the remainder of this paper, the transpose of the exact flux Jacobian matrix will be denoted $A = \left[ \frac{ \partial \mathbf{R}}{ \partial \mathbf{W}} \right]^T$.

\section{Minimal residual Krylov subspace methods combined with spectral deflation}

In this section we focus on a particular minimal residual norm Krylov subspace method for the solution of linear systems with a non-symmetric real coefficient matrix of type 
\begin{equation}
    \label{eq:Linear_System}
    Ax = b, \qquad A \in \mathbb{R}^{N \times N}; \quad b, \hspace{0.1cm} x \hspace{0.1cm} \in \mathbb{R}^N
\end{equation}

A right-preconditioned system is considered so that system (\ref{eq:Linear_System}) becomes
\begin{align}
    \label{eq:Varaible_Precond_System}
    A\mathcal{M}(t) = b, \\
    \label{eq6}
    x = \mathcal{M}(t)
\end{align}
with $t \in \mathbb{R}^{N}$ and $\mathcal{M}:\mathbb{R}^{N} \rightarrow \mathbb{R}^{N}$ the preconditioning operator which may be a nonlinear function.

\subsection{Flexible GMRES algorithm with right preconditioning}

Saad has proposed a minimal residual norm subspace method based on the standard GMRES approach \cite{saad1986} that allows a variable nonlinear preconditioning function $\mathcal{M}_{j}:\mathbb{R}^{N} \rightarrow \mathbb{R}^{N}$ at each iteration $j$, \cite{saad1993flexible}.

Starting from an initial guess $x_{0}$, the flexible Arnoldi relation is written as:  
\begin{equation}
    \label{eq:Flexible_Arnoldi}
    AZ_{m} = V_{m+1}\Bar{H}_{m},
\end{equation}

\noindent where the matrices $V_{m+1} \in \mathbb{R}^{N \times (m+1)}, Z_{m} \in \mathbb{R}^{N \times m}$ and $\bar{H}_{m} \in \mathbb{R}^{(m+1) \times m}$ stand for the orthonormal basis of the Krylov space, the solution space and the upper Hessenberg matrix respectively. The approximate solution is written as $x_{m} = x_{0} + Z_{m}y_{m}$ where $y_{m}$ minimizes $|| r_{0} -AZ_{m}y ||_{2}$ over $x_{0} + span\lbrace Z_{m} \rbrace$, with $Z_{m}=\mathcal{M}V_{m}=\left[\mathcal{M}_{1}(v_{1}),...,\mathcal{M}_{m}(v_{m}) \right]$ and both $Z_{m}$ and $V_{m}$ need to be stored. We point out that the operator $\mathcal{M}$ represents the action of the nonlinear operators $\mathcal{M}_{j}$ on the set of basis vectors $v_{j}$.

The restarted FGMRES($m$,$m_{i}$) algorithm is presented in Algorithm \ref{alg:FGMRES}. We denote by $m_{i}$ the size of the Krylov subspace associated to the GMRES solver devoted to the inner linear system. We point out that the stopping criterion is essentially based on the true relative residual defined by $\frac{||b-Ax||}{||b||}$ (see \cite{fraysse2008algorithm} for several definitions of stopping criteria and some practical considerations for the implementation of the key points of the algorithm). Still from an implementation point of view, the relative least-squares residual is currently used as a cheap approximation of this quantity at each Arnoldi iteration.

\begin{algorithm}[H]
\caption{FGMRES($m$,$m_{i}$)}
\begin{algorithmic}[1]
\STATE  Choose an initial guess $x_{0}$ and a convergence threshold $\epsilon$
\STATE  Compute $r_{0} = b-Ax_{0}$, $\beta =||r_{0}||$ and  $v_{1} = r_{0}/\beta$ ;
\STATE  c = $[\beta,0...0]^{T} \in \mathbb{R}^{m+1}$
\FOR{$j=1,...,m$}
    \STATE $z_{j} = \mathcal{M}_{j}(v_{j})$
    \STATE $w = Az_{j}$
    \FOR{$i=1,...j$}
        \STATE $h_{i,j} = v_{i}^{T}w$
        \STATE $w = w -h_{i,j}v_{i}$
    \ENDFOR
    \STATE $h_{j+1,j} =||w||$ and $v_{j+1} = w/h_{j+1,j}$
    \STATE Solve the least-squares problem $min_{y} ||c - \bar{H_{j}}y||$ for $y^{*}$
    \STATE Exit if $||c - \bar{H}_{m}y^{*}||/||b|| \le \epsilon $
\ENDFOR
\STATE Compute $x_{m} = x_{0} + Z_{m}y^{*}$ where $Z_{m} = [z_{1},...,z_{m}]$
\STATE Set $x_{0} = x_{m}$ and go to 2
\end{algorithmic}
\label{alg:FGMRES}
\end{algorithm}

The relative true residual is only computed at the end of each cycle and is used to construct the first vector of the next Krylov subspace basis.
In the case of large and ill-conditioned linear systems, least-squares and true residuals may differ due to loss of orthogonality during the construction of the Krylov basis.
A standard way to tackle such a phenomenon is to ask for a second iteration of the Modified Gram-Schmidt algorithm (loop from line 7 to 10) in order to reinforce the orthogonality of the Krylov basis.
In this paper, the nested GMRES strategy is adopted for the numerical experiments. Two preconditioning strategies are considered for the inner GMRES. The first one consists in a block version of a standard LU-SGS iterative solver. LU-SGS is applied to a first order diagonally dominant upwind approximation of the flux Jacobian matrix inspired by~\cite{yoon1988}. This operator is based on a first order spatial discretization of the convective and of the viscous fluxes using a simplifying thin layer assumption~\cite{peter2007large}. This strategy leads to a very compact stencil for the preconditioning matrix which will be denoted by $\textbf{J}^{APP}_{O1}$ in the sequel of the paper. The second one is a Block Incomplete LU (BILU($k$)) factorization applied to either an approximate or exact flux Jacobian matrix. For the so-called first-order exact Jacobian matrix $\textbf{J}^{EX}_{O1}$ a first-order spatial Roe scheme is used for the discretization of the mean-flow convective fluxes and a 5-point corrected centered discretization scheme is used for the diffusive fluxes.  More specifically, the spatial gradients at the cell interfaces are modified to avoid high frequency oscillations (see~\cite{resmini2015sensitivity} or~\cite{puigt2014}). The BILU($k$) preconditioner will be applied either to the first-order approximate Jacobian matrix $\textbf{J}^{APP}_{O1}$, or to the first-order exact Jacobian matrix $\textbf{J}^{EX}_{O1}$. $\textbf{J}^{EX}_{O1}$ is different from $\textbf{J}^{APP}_{O1}$ when it comes to memory footprint. More specifically, $\textbf{J}^{EX}_{O1}$ has a 9-point stencil in 2D whereas a 5-point stencil is associated with $\textbf{J}^{APP}_{O1}$. In 3D, we have a 7-point stencil for $\textbf{J}^{APP}_{O1}$ and a stencil of 19 points for $\textbf{J}^{EX}_{O1}$. Consequently, a better robustness is achieved but at the price of about twice the storage for $\textbf{J}^{EX}_{O1}$ compared to $\textbf{J}^{APP}_{O1}$. For the high-order DG formulation, an exact third-order Jacobian matrix $\textbf{J}^{EX}_{O3}$ is built. Only the BILU(0) preconditioner is considered in the latter case. The size of the matrix reads $N_{elt}N_{eq}N_{p}$ where $N_{elt}$, $N_{eq}$, and $N_{p}$ represent the numbers of mesh elements, of equations, and of degrees of freedom per equation respectively. We have $N_{p}=1$ for the FV computations and either $N_{p}=6$ in 2D, or $N_{p}=10$ in 3D for the DG computations. Ignoring boundary conditions, the number of nonzero entries in the matrix can be estimated through $NNZ=(1+2sd)N_{elt}(N_{eq}N_{p})^{2}$ where $1+2sd$ is the number of points in the stencil and $d$ the space dimension. For the DG computations, we always have $s=1$, but $N_p>1$ which usually results in larger $NNZ$ values compared to FV even when $N_{elt}$ is lower. Both matrices are sparse and contain dense blocks of size $N_{eq}N_{p}$ which will have a strong impact on the performances of the linear solvers as shown in the numerical experiments. 

We point out that the relevant numerical ingredients that characterize the GMRES algorithm are the matrix product (step 5 in Algorithm \ref{alg:FGMRES}), the preconditioning step (step 4) and the scalar product (step 7). These algebraic operations are global in conjunction with a domain decomposition method. More specifically, the globalization of the preconditioner (step 4) is achieved with a Restricted Additive Schwarz method~\cite{cai1999}. In addition, the product by the operator $A$ (step 5) is exact. We thus get a global and parallel FGMRES($m$,$m_{i}$). 

\subsection{Domain decomposition-based preconditioners}
We now describe two stationary block preconditioners that are used in this work and that have been enhanced with overlapping techniques. The first one is a combination of a Restrictive Additive Schwarz (RAS) domain decomposition method and of a LU-SGS relaxation strategy. It will be applied only to the FV problem. The second one a is a BILU(0) factorization over and extended partition of the fluid domain by some level of overlapping and will be applied to the DG problem.

\subsubsection{Overlap Restrictive Additive Schwarz LU-SGS preconditioner}
\label{lussor}
The LU-SSGS solver is used as a stationary preconditioner for the structured multi-block FV adjoint solution.
It consists in an iterative technique combining a LU (Lower-Upper) factorization with a relaxation method~\citep{saad2003}. In practice, we set the LU-SSOR relaxation parameter $\omega$ to $1$, and the method is then referred to as LU-SGS. In the context of preconditioning, we are interested in the approximate solution of the linear system $Mx=v$. In the factorization step, the implicit operator $M$ is decomposed into a lower triangular matrix $L$, a diagonal matrix $D$ and an upper triangular matrix $U$: $M=L+D+U$.


The solution of the system is then approximated by relaxation~\citep{peter2007large,yoon1988}, performing forward and backward successive sweeps.
Each relaxation cycle involves two stages:
\begin{align}
(U+D)x^{p+\frac{1}{2}} &= v -Lx^{p} \label{eq_LUSGS3} \\
(L+D)x^{p+1} &= v -Ux^{p+\frac{1}{2}} \label{eq_LUSGS4}
\end{align}
where $p\in\mathbb{N}$ is index of the current cycle, starting with $x^0=0$.
For one cycle of relaxation we can explicit the corresponding preconditioning operator as $M_{SGS}=(U+D)D^{-1}(L+D)$.
It is worth mentioning that in the preconditioning context, typically $2$ to $4$ cycles (of 2 sweeps) are enough to provide a useful approximate solution.
Besides, as we deal with block-structured problems, a block LU-SGS algorithm has been implemented, i.e., $D$ is block-diagonal. As already mentioned, the matrix $\textbf{J}^{APP}_{O1}$ is diagonally dominant by construction. However, when the turbulence model is linearized our numerical experiments have shown that it is beneficial to add a diagonal of the form $I/\Delta \tau$ where $\Delta\tau$ is consistent with a standard local pseudo-time step typically introduced in the implicit stage of a nonlinear inexact Newton steady solver. The local pseudo-time step is computed from a prescribed target Courant-Friedrichs-Lewy (CFL) number and from the steady-state flow field. This requires a trial and error strategy to choose the right CFL number but, as we will see, an appropriate CFL can dramatically accelerate the convergence of the adjoint solver.

In a parallel multi-domain coupling context, the part of the right-hand side of (\ref{eq_LUSGS3},~\ref{eq_LUSGS4}) associated with overlapping cells is updated before each stage in order to account for the information from connected blocks. For the current stage, depending on the ordering of the blocks, the right-hand side is built from the solution at the previous stage or from the blocks whose solution is already up-to-date. As we deal with first order flux Jacobian matrices, the size of the overlap is one. In the remainder of this paper the notation LUSGS(N,CFL) refers to N stages of sweep, that is, N is an even number.

\subsubsection{Improved overlap BILU(0) preconditioner}

For the unstructured DG adjoint solution, the BILU(0) preconditioner used in the inner GMRES solver is combined with the RAS method.
Such a strategy aims to enhance the convergence of the inner system when the number of partitions increase and therefore the overall convergence of the outer system.
Interactions between neighboring partitions are taken into account in order to enlarge the slice of the global implicit matrix that is carried out per partition.
As the DG formulation requires a compact stencil of size one, we only need to exchange blocks from our direct neighbors by MPI.
The BILU(0) algorithm is then applied on this slice to construct an enlarged local preconditioning matrix.
To apply the preconditioning step on a global vector, a restriction operator is first used to obtain the vector corresponding to the enlarged slice of the current partition thanks to an extra MPI communication. At the end, no prolongation operator is played and we only retain entries of the vector associated with the original slice of the current partition.
The matrix-vector product with the global implicit matrix is still performed in parallel in an exact way.

\subsection{Deflation approach}

The main drawback of the restarted GMRES($m$) is the loss of spectral information contained into the current Krylov subspace during the restarting procedure.
Let us recall the definition of a Ritz pair \citep{morgan1998harmonic} as it plays an important role in the strategy of deflated restarting.
The standard Rayleigh-Ritz problem specifies an orthogonality condition on the spectral residual:

\begin{equation}
(Ay -\lambda y) \perp \mathcal{K}_{m}(A,r_{0}) \quad \forall y \in \mathcal{K}_{m}(A,r_{0})
\label{eq:Ritz_problem}
\end{equation} 

Using $y=V_{m}g$ and the standard Arnoldi relation $AV_{m} = V_{m+1}\bar{H}_{m}$ in (\ref{eq:Ritz_problem}), we have 

\begin{equation}
H_{m}g = \lambda g
\label{eq:Ritz_spectral}
\end{equation}

Thus, the spectral residual norm in (\ref{eq:Ritz_problem}) of the Ritz pair $\left\{\lambda,y=V_{m}g\right\}$ satisfies:
\begin{equation}
\label{eq:Residual_eigenvalue_norm}
  \|A(V_{m}g) - \lambda (V_{m}g)\| = \|A(V_{m}g)-V_{m}H_{m}g\| = \|V_{m+1}\Bar{H}_{m}g - V_{m}H_{m}g\| = | h_{m+1,m} | |e_{m}^{T}g|
\end{equation}
From a small value of $|h_{m+1,m}| |e_{m}^{T}g|$, the Arnoldi method takes as an approximate eigenvalue-eigenvector pair of the operator $A$ the Ritz pair~\cite{arnoldi1951,saad1980}.
Indeed, neglecting the last row of the rectangular upper Hessenberg matrix $\Bar{H}_{m}$ leads to:

\begin{equation}
    \label{eq:Hessenberg_Restriction_operator}
    H_{m} = V_{m}^{T}AV_{m}
\end{equation}

Therefore, the spectrum of $H_{m}$ naturally approximates a part of the spectrum of $A\mathcal{M}$.
The idea of deflation techniques is to keep relevant spectral information from $H_{m}$ in the search space of the next cycle to expect a better convergence of the Krylov iterative methods.
In~\cite{giraud2010}, Giraud et al. take into account both smallest and largest eigenvalues to maximize the deflation effect.
In contrast, Morgan~\cite{morgan2002} only deflates the smallest ones. This last strategy will be adopted for our numerical experiments.


Actually, Ritz values of the operator $A$ give a good approximation of its exterior eigenvalues.
Unfortunately, interior eigenvalues are of greater interest because they are generally responsible for the convergence stagnation.
The harmonic Ritz values of $A$ are defined as the Ritz values of $A^{-1}$ with respect to the subspace $A\mathcal{K}_{m}(A,r_{0})$.
The resulting orthogonality condition can be expressed as:
\begin{equation}
	\label{eq;Harmonic_Eigenvalue_equation}
    (A^{-1}y - \theta y) \perp A\mathcal{K}_{m}(A,r_{0}) \quad \forall y \in A\mathcal{K}_{m}(A,r_{0})
\end{equation}
The orthogonality condition for harmonic Ritz values leads to the generalized eigenvalue problem:
\begin{equation}
	\label{eq:Generalized_eigenvalue_problem}
    \theta \Bar{H}_{m}^{T} \Bar{H}_{m}g = H_{m}^{T} g,
\end{equation}
After some algebraic manipulations, (\ref{eq:Generalized_eigenvalue_problem}) can be reformulated as a standard eigenvalue problem:
\begin{equation}
	\label{eq:Standard_eigenvalue_problem}
    (H_{m} + h_{m+1,m}^{2}H_{m}^{-T}e_{m}e_{m}^{T})g = \theta^{-1} g.
\end{equation}
where $\theta^{-1}$ is the harmonic Ritz value and the corresponding harmonic Ritz vector is $y=AV_{m}g$.
As harmonic Ritz values of $A$ are Ritz values of $A^{-1}$, exterior eigenvalues $\theta$ of $A^{-1}$ are now supposed to be well approximated.
Therefore, theirs inverses $\theta^{-1}$ are expected to be good approximations of eigenvalues of $A$ in the neighborhood of zero.

Obviously, in exact arithmetic solutions to (\ref{eq:Generalized_eigenvalue_problem}) and (\ref{eq:Standard_eigenvalue_problem}) are identical. But it is not the case in finite precision since the operator $\Bar{H}_{m}^{T} \Bar{H}_{m}$ is usually ill-conditioned in the fully linearized turbulence case. Therefore, the accurate estimation of the eigenvectors could be strongly altered and lead to stagnation of the relative true residual of the GMRES process. In contrast, we get better results at the conditioning level with the standard eigenvalues problem overcoming the stagnation of the relative true residual.

As mentionned by Giraud et al. \cite{giraud2010}, the flexible Arnoldi relation obtained at each restart within the FGMRES with deflated restarting (FGMRES-DR) framework given by 

\begin{equation}
\label{eq:Flexible_Arnoldi_Realtion_Deflation}
AZ_{k} = V_{k+1}\bar{H}_{k},
\end{equation}

\noindent holds with $Z_{k} = Z_{m}P_{k}$, $V_{k+1} = V_{m+1}P_{k+1}$ and $\bar{H}_{k} = P_{k+1}\bar{H}_{m}P_{k}$ where $P_{k} \in \mathbb{R}^{m \times k}$ corresponds to the orthonormal matrix whose columns are spanned by the eigenvectors of (\ref{eq:Standard_eigenvalue_problem}).

Also, the right-hand side of the least-squares problem is computed at each restart as $c = V_{m+1}^{T}r_{m}$ which requires $2N(m+1)$ operations. Rollins and Fichtner \citep{rollin2008improving} have proposed an efficient way to compute $c$ so that we can save some inner products. Indeed, they demonstrate that the residual $r_{m}$ is a linear combination of the columns of the deflation subspace $V_{k+1}$. Consequently, we have $c=V^{T}_{k+1}r_{m}$ with $2N(k+1)$ operations. Also, they improved the construction of $P_{k}$ and in particular the last column of this matrix which is usually chosen as the vector $c - \bar{H}_{m}y_{m}$. More specifically, they demonstrated that the vector $c - \bar{H}_{m}y_{m}$ is colinear to the vector $[-\beta f^{T} \; 1]^{T}$ with $f=H^{-T}_{m}e_{m}$. The explicit relation is given below:

\begin{equation}
\label{eq:Rollin_vector}
c - \bar{H}_{m}y_{m} = \begin{pmatrix} -\beta f \\ 1 \end{pmatrix} \left( \frac{\omega - \beta f^{T}v}{1+\beta^{2}f^{T}f} \right),
\end{equation}

\noindent where $v$ and $\omega$ are respectively the first $m$ rows of $c$ and the last element of $c$.

In exact arithmetic, using the vector $[\beta f^{T} \; 1]^{T}$ or the vector $c -\bar{H}_{m}y_{m}$ as the last column of $P_{k}$ is equivalent due to the colinearity property. However, in finite precision arithmetic, the former vector is numerically preferable. More precisely, the use of the vector $f$ both to compute the eigenvectors of $H_{m} + \beta^{2}fe^{T}_{m}$ and to construct the last column of matrix $P_{k}$ enables to reduce the rounding errors in the Arnoldi relation after a restart (see \citep{rollin2008improving}). Consequently, maintaining the Arnoldi relation to a high accuracy leads to a better convergence of the GMRES solver. All these numerical considerations have been illustrated in the next section.

\section{Numerical experiments}

To evaluate this GMRES solver capabilities, two dedicated test cases have been defined. The first test case was conducted on an adjoint system of equations based on the turbulent transonic RANS flow over a two-dimensional ONERA OAT15A airfoil. The Mach number is 0.734, the Reynolds number based on the chord and freestream conditions is 6.5 $\times$ $10^{6}$ and the angle of attack is 1.15°. The second test case is based on an adjoint system of equations for the three-dimensional turbulent transonic RANS flow over the ONERA M6 wing at Mach number $M=0.84$, at Reynolds number, based on the mean chord and freestream conditions, of 1.17 $\times$ $10^{7}$ and with an angle of attack of 3.06°. In both computations the Spalart Allmaras one equation turbulence model has been selected. We focus on the solution of the adjoint system with the fully linearized turbulence model. For the FV solver, the numerical scheme is a second order upwind Roe spatial discretization associated with a MUSCL reconstruction and a van Albada limiter. For the preconditioning step $\textbf{J}^{EX}_{O1}$ and $\textbf{J}^{APP}_{O1}$ are considered. For the DG solver, the numerical scheme relies on the Roe flux for the convective term and on the Bassi and co-workers (BR2)~\cite{bassi1997high} formulation for the diffusive term, and $\textbf{J}^{EX}_{O3}$ is used as the preconditioner.

We emphasize that the FGMRES solvers have been implemented in different computing environments. For instance, the FGMRES solver dedicated to the FV problem has been developed in a modular Python framework. The characteristic algebraic operations (i.e scalar product, preconditioning step and matrix product) have been directly used by wrapping, via the open source software tool SWIG, the corresponding routines from the kernel of the CFD solver $\textit{elsA}$. This organization is very flexible as it allows a rapid prototyping of promising numerical strategies in a Python environment while still achieving a good level of performance. In contrast, the FGMRES solver dedicated to the DG problem has been implemented in a Fortran environment with a different data structure. Therefore, no comparison can be drawn from a high-performance computing point of view. 

\subsection{ONERA OAT15A airfoil test case}

The ONERA OAT15A airfoil has a thick trailing edge such that the initial structured mesh is made of two blocks. The structured mesh for the FV calculation is composed of 158 208 quadrilaterals and the unstructured mesh for the DG one of 64 416 quadrilaterals.

\begin{figure}[H]
    \centering
    \begin{subfigure}[b]{.5\linewidth}
    \centering
    \includegraphics[width=.8\linewidth]{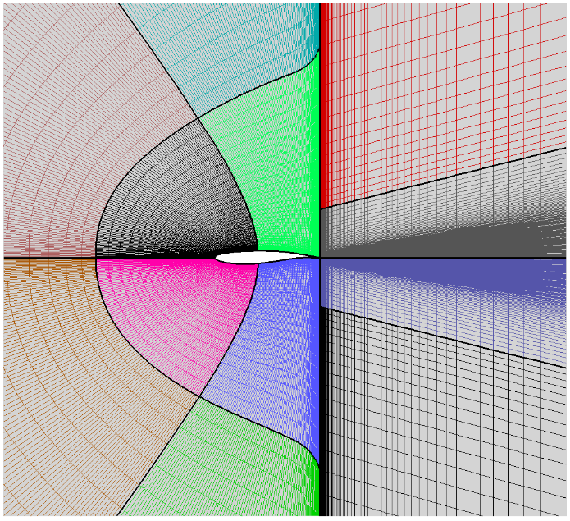}
    \caption{}\label{fig1:a}
    \end{subfigure}
    \begin{subfigure}[b]{.5\linewidth}
    \centering
    \includegraphics[width=.8\linewidth]{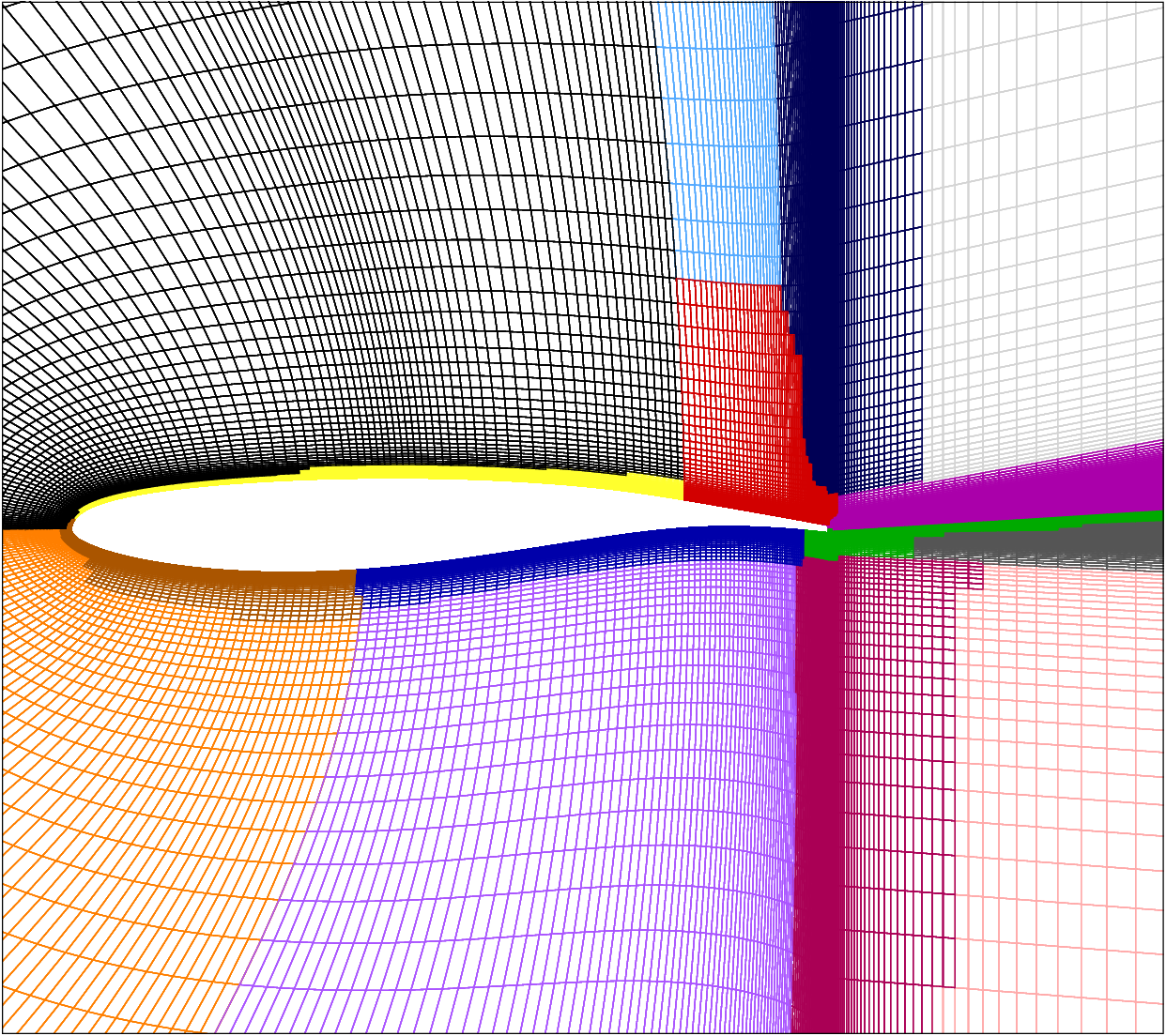}
    \caption{}\label{fig1:b}
    \end{subfigure}
    \begin{subfigure}[b]{.5\linewidth}
    \centering
    \includegraphics[width=.8\linewidth]{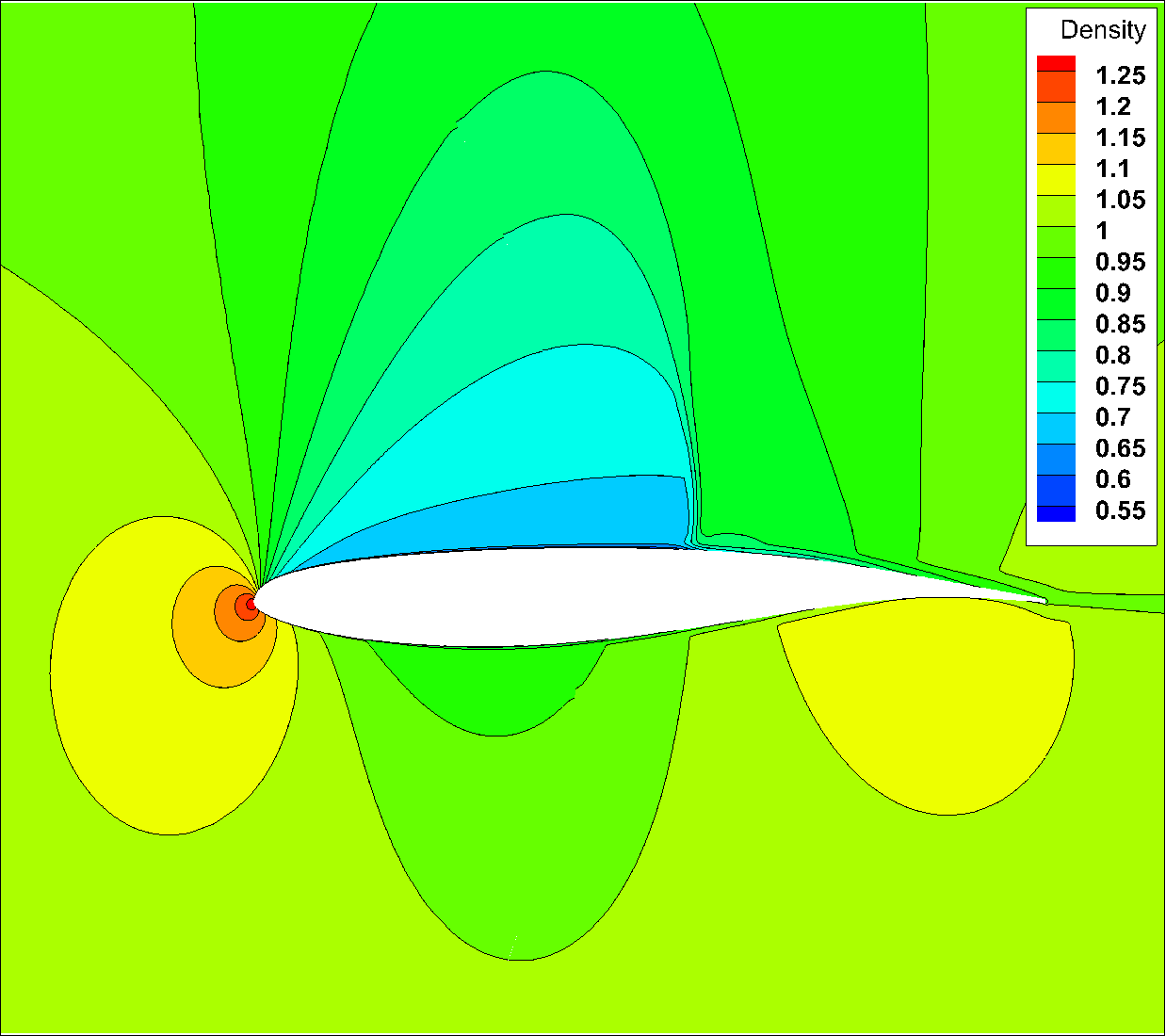}
    \caption{}\label{fig1:c}
    \end{subfigure}
    \begin{subfigure}[b]{.5\linewidth}
    \centering
    \includegraphics[width=.8\linewidth]{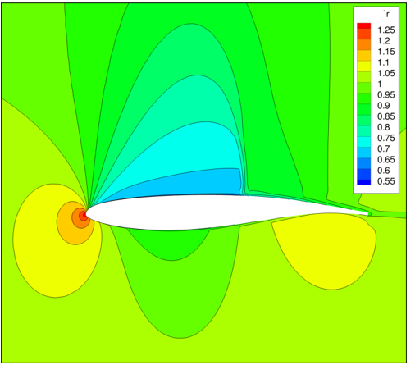}
    \caption{}\label{fig1:d}
    \end{subfigure}
    \caption{2D transonic turbulent flow over the OAT15A airfoil: (a) and (b) illustrate the structured and unstructured meshes partitioned into 16 subdomains. The corresponding steady-state density field is illustrated in (c) and (d) for the FV and DG schemes respectively.}\label{fig1}
\end{figure}

Figure~\ref{fig1:a} and Figure~\ref{fig1:b} depict the mesh partitioning obtained on 16 cores for the scalability study.

For the unstructured case, a general matrix-based algorithm~\cite{mavri2020} is implemented to construct lines by grouping strongly-connected degrees of freedom.
The definition of these lines relies on a directional stiffness measure extracted from a stiffness matrix that represents the strength of couplings between the unknowns.
Then, a multilevel k-way graph partitioning
algorithm from the METIS library~\cite{metis1999} is applied based on a strategy taking into account weights on both the vertices and the edges of the graph. Finally, four weights have been defined and calibrated to distinguish mesh cells belonging to a strong connection line but also mesh edges whose one of the parent cell is affected. The main interest is to preserve as much as possible certain critical areas of the original domain from being arbitrarily partitioned. For example in the case of a shock wave propagation or a turbulent boundary layer development. The resulting largest load imbalance for the unstructured case finally is of the same order as that obtained for the structured case and is about 8 $\%$.

\subsubsection{Adjoint sensitivity analysis}
\label{sec4_2}
Recall that we want to solve the adjoint system (\ref{eq:Adjoint_System}) efficiently by using the inner-outer GMRES solver with deflated restarting denoted by FGMRES-DR($m$,$m_{i}$,$k$). We adopt the relevant numerical parameters in Table~\ref{table1}. The inner GMRES solver will not be restarted as only small sizes of the Krylov inner subspace will be considered.

\begin{table}[ht]
    \centering
    \small
    \begin{tabular}{|c|c|c|}
        \hline
         $m$ & 60 & Size of the outer Krylov space \\
        \hline
         $k$ & 20 & Number of deflated vectors \\
         \hline
         tol outer& 1.e-9 & Relative convergence threshold \\
         \hline
         $m_{i}$ & 20 & Maximum size of the inner Krylov space \\
         \hline
         tol inner & 0.5 & Convergence threshold for the inner GMRES \\
         \hline
    \end{tabular}
    \caption{FGMRES-DR relevant numerical parameters}\label{table1}
\end{table}

\begin{figure}[H]
    \centering
    \begin{subfigure}[b]{.5\linewidth}
    \centering
    \includegraphics[trim=0.25cm 0.5cm 1.5cm 1.5cm,clip,width=1.05\linewidth]{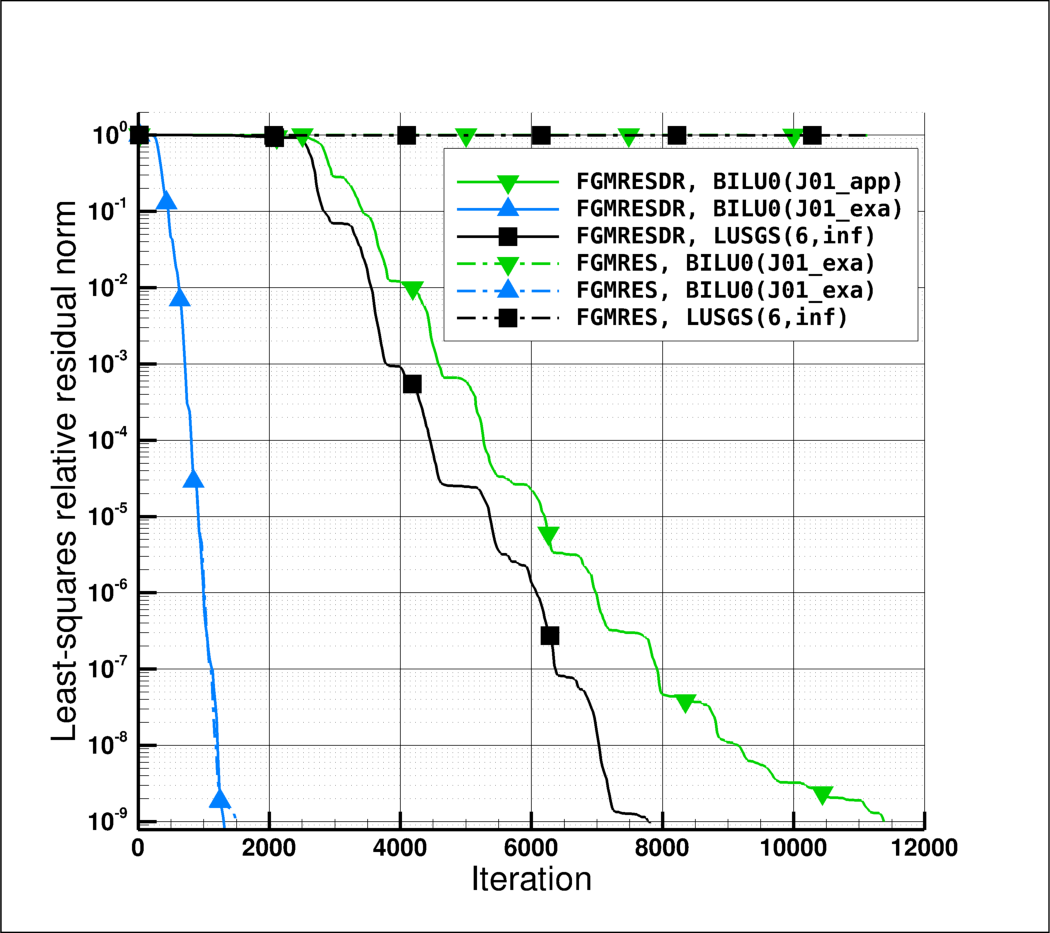}
    \caption{FV case}\label{fig2:a}
    \end{subfigure}
    \begin{subfigure}[b]{.5\linewidth}
    \centering
    \includegraphics[trim=0.25cm 0.5cm 1.5cm 1.5cm,clip,width=1.05\linewidth]{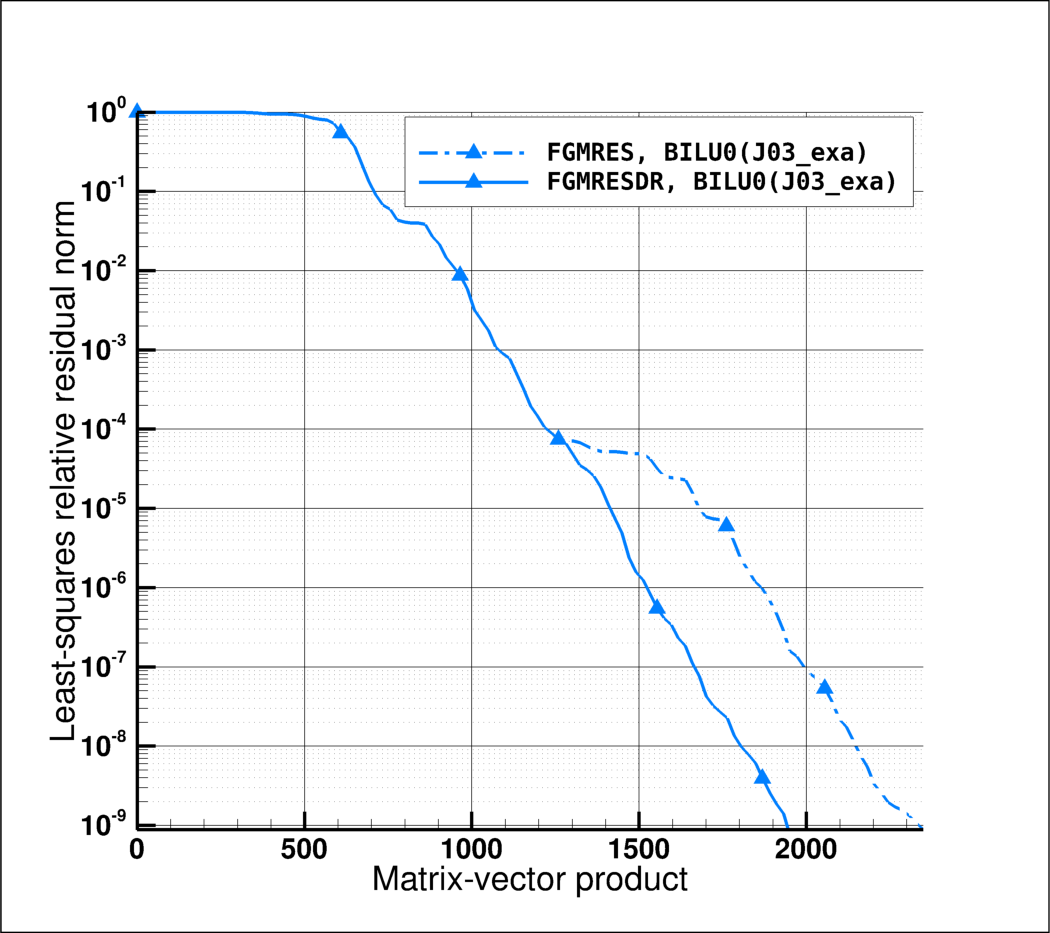}
    \caption{DG case}\label{fig2:b}
    \end{subfigure}
    \caption{Impact of deflation for various preconditioners. The relative residual norm convergence history is plotted with respect to iterations. The impact of deflation is dramatically effective both on the LUSGS and BILU preconditioners applied to the first-order approximate Jacobian matrix. We recall the numerical parameters of the FGMRES-DR solver as follows: $m=60, m_{i}=20$ and $k=20$.}\label{fig2}
\end{figure}

All computations are performed with 16 cores on 16-blocks structured and unstructured grid thus with one domain associated with each core. The underlying algebraic problem size is fixed to 1 million for the FV scheme and 2 million for the DG scheme. The number of entries is about 35 million and 288 million for the FV scheme and the DG scheme respectively. Figure~\ref{fig2} displays the relative least-squares residual norm convergence histories of FGMRES-DR for both FV and DG schemes. In the FV case, we notice a stagnation of FGMRES(60,20) for both LU-SGS and BILU(0) algorithms applied to the first-order approximate Jacobian matrix $\textbf{J}^{APP}_{O1}$. To overcome this stagnation it is natural to extend the approximation space $\mathcal{K}_{m}$ in order to capture as much as possible the eigenvalues of $A\mathcal{M}$. Since we cannot afford to extend the projection space due to memory limitation, the deflation strategy remains the best alternative to retrieve the residual convergence of a quasi-full FGMRES solver at low computational cost. A deflation of 30$\%$ of $m$ is sufficient and the resulting FGMRES-DR(60,20,20) solver, preconditioned by LU-SGS(6) or BILU(0) now converges as can be seen in Figure~\ref{fig2}a. However, when it comes to use the first order exact Jacobian matrix $\textbf{J}^{EX}_{O1}$, FGMRES(60,20) already converges in a few iterations. In this case, deflating does not lead to any improvements which illustrates how robust $\textbf{J}^{EX}_{O1}$ is. The downside of this is that we need approximately twice the storage of $\textbf{J}^{APP}_{O1}$. Figure~\ref{fig2}b shows convergence histories of both FGMRES(60,20,20) and FGMRES-DR(60,20,20) for the DG case. From a DG formulation, resulting Jacobian matrices are often ill-conditioned, especially when the spatial order of the scheme becomes high. A variant of the BILU(0) factorization applied to expanded domains (by some level of overlapping) has been implemented in order to globalize as much as possible the preconditioning effect on the whole fluid domain. In both cases, the residuals converge in a few iterations and deflation is shown to slightly improve the convergence rate. We have reported the performance details of the two solvers in Table~\ref{table2} and Table~\ref{table3} for both FV and DG problems.

The number of nonzero entries (NNZ) of the preconditioner is reported in Table~\ref{table2} and Table~\ref{table3}.

\begin{table}[H]
    \centering
    \begin{tabular}{|c|c|c|c|c|} 
    \cline{1-5}
    \multicolumn{1}{|c|}{$\mathcal{M}$ = BILU(0)} & \multicolumn{2}{|c|}{FGMRES(60,20)} & \multicolumn{2}{|c|}{FGMRES-DR(60,20,20)} \\
    \cline{1-5}
    \multicolumn{1}{|c|}{Jacobian matrix} & $\textbf{J}^{EX}_{O1}$ &  $\textbf{J}^{APP}_{O1}$ &  $\textbf{J}^{EX}_{O1}$ & $\textbf{J}^{APP}_{O1}$ \\ 
    \hline
    NNZ (million) & 28 & 21 & 28 & 21 \\
    \hline 
    $\#$ its & 74 & - & 74 & 546 \\
    \hline 
    $\#$ Mvps & 1498 &  - & 1321 & 11381 \\
    \hline
    \end{tabular} 
    \caption{Performance of FGMRES-DR for the FV method}\label{table2}
\end{table}

\begin{table}[ht]
\centering
    \begin{tabular}{|c|c|c|} 
    \cline{1-3}
    \multicolumn{1}{|c|}{$\mathcal{M}$ = BILU(0)+RAS} & \multicolumn{1}{|c|}{FGMRES(60,20)} & \multicolumn{1}{|c|}{FGMRES-DR(60,20,20)} \\
    \cline{1-3}
    \multicolumn{1}{|c|}{Jacobian matrix} &  $\textbf{J}^{EX}_{O3}$ &  $\textbf{J}^{EX}_{O3}$ \\ 
    \hline
    NNZ (million) & 288 & 288 \\
    \hline 
    $\#$ its & 111 & 93 \\
    \hline 
    $\#$ Mvps & 2500 &  1940 \\
    \hline
    \end{tabular}
    \caption{Performance of FGMRES-DR for the DG method}\label{table3}
\end{table}

\subsubsection{Scalability}

We have demonstrated in the previous section how competitive the deflated inner-outer GMRES is when applied to stiff problems. As already mentioned, the inner GMRES is by construction a global preconditioner and we are then interested in how its performances scale in a parallel environment. Here, we pay a special attention to the strong scalability in the sense that the problem size is fixed while varying the number of cores. 

\begin{figure}[ht]
    \centering
    \begin{subfigure}[b]{.5\linewidth}
    \centering
    \includegraphics[trim=0.25cm 0.5cm 1.5cm 1.5cm,clip,width=1.05\linewidth]{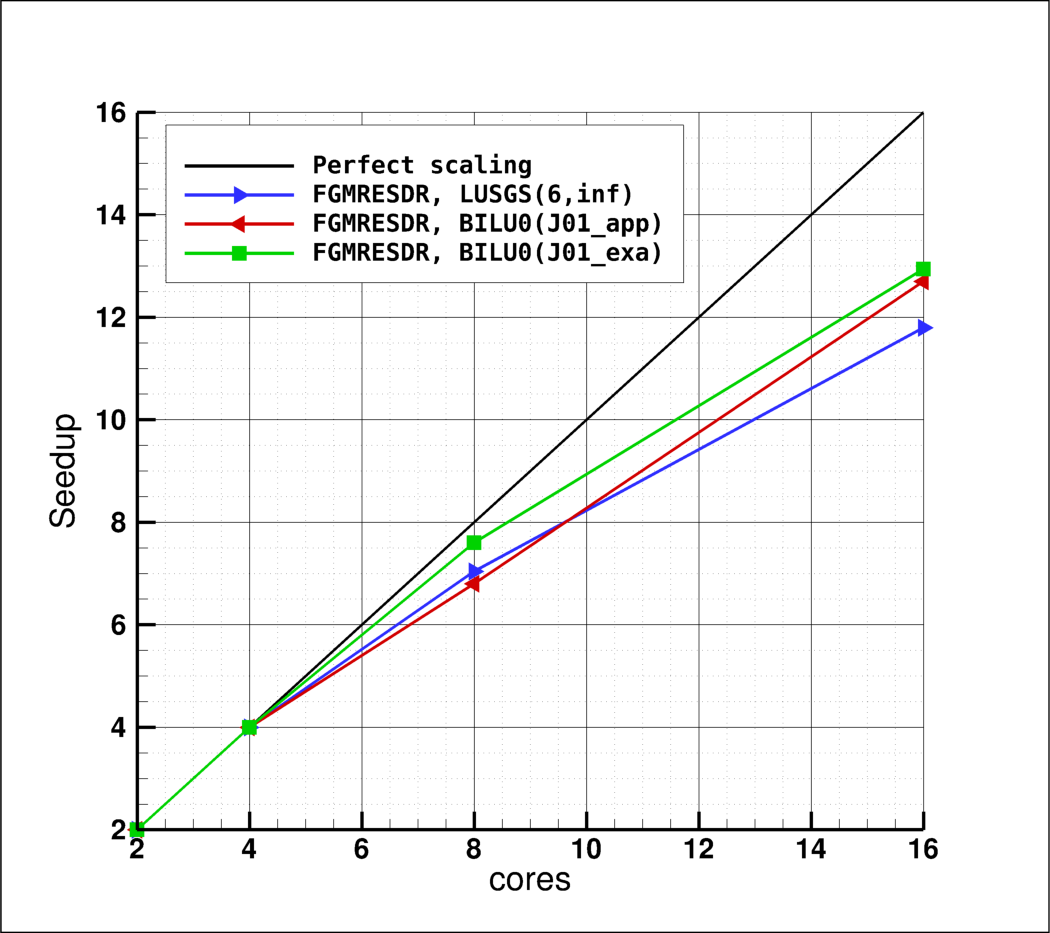}
    \caption{FV scheme}\label{fig3:a}
    \end{subfigure}
    \begin{subfigure}[b]{.5\linewidth}
    \centering
    \includegraphics[trim=0.25cm 0.5cm 1.5cm 1.5cm,clip,width=1.05\linewidth]{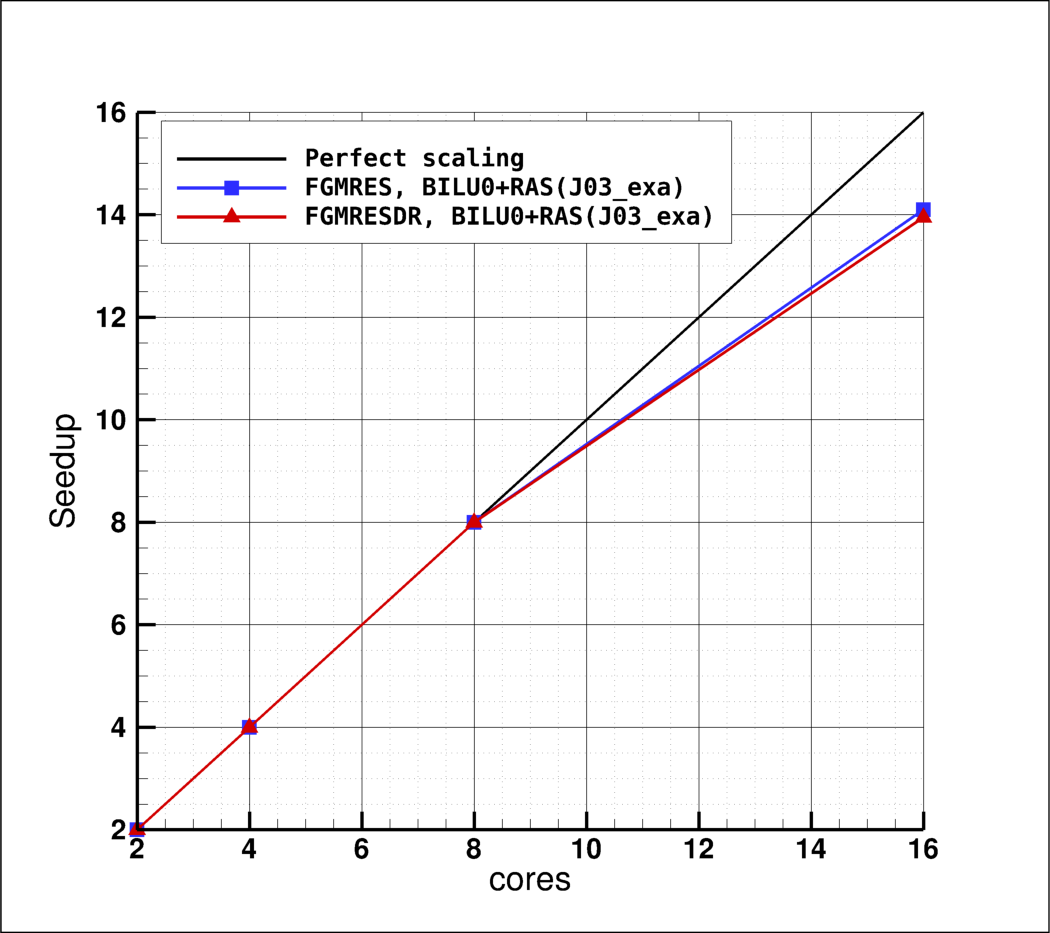}
    \caption{DG scheme}\label{fig3:b}
    \end{subfigure}
    \caption{Strong scalability analysis for FGMRES-DR(60,20,20) solvers on the OAT15A adjoint system.}\label{fig3}
\end{figure}

For the FV scheme, we notice in Figure~\ref{fig3:a} that the efficiency reaches up to 75 $\%$ for the LU-SGS(6) applied to $\textbf{J}^{APP}_{O1}$ while the speedup for the BILU(0) applied to $\textbf{J}^{EX}_{O1}$ reaches up to 80 $\%$. The same results have been observed for the DG scheme with BILU(0) and RAS algorithms since an efficiency of 85 $\%$ for FGMRES(60,20) as well as for the FGMRES-DR(60,20,20) is measured (Figure~\ref{fig3:b}). We point out that the local nature of the BILU(0) preconditionner does not have a real impact on the robustness of the inner-outer GMRES methods unlike the standard GMRES~\cite{simoncini2002flexible}. Indeed, one of the main advantages of the inner-outer GMRES is the global nature of the inner GMRES preconditionner since it solves the same initial linear system (4).

For both structured and unstructured cases, the number of iterations to converge the adjoint problem shows a variation up to 20$\%$ compared to the two-subdomain reference configuration.

\subsection{ONERA M6 wing test case}

\begin{figure}[ht]
    \centering
    \begin{subfigure}[b]{.5\linewidth}
    \centering
    \includegraphics[width=.8\linewidth]{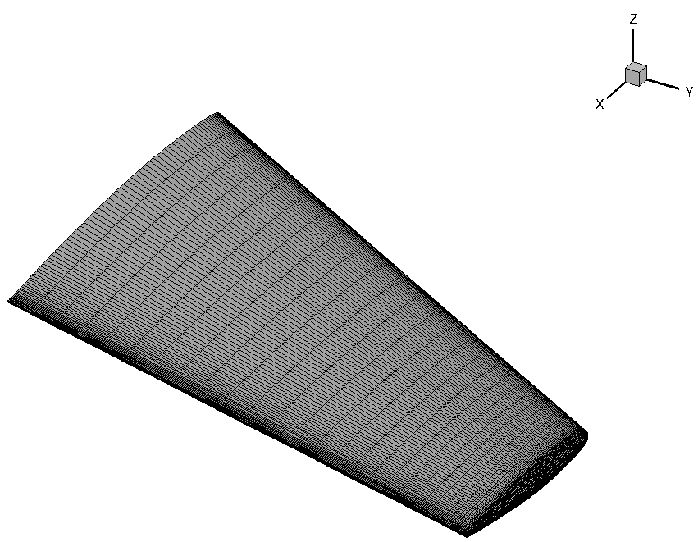}
    \caption{}\label{fig4:a}
    \end{subfigure}
    \begin{subfigure}[b]{.5\linewidth}
    \centering
    \includegraphics[width=.8\linewidth]{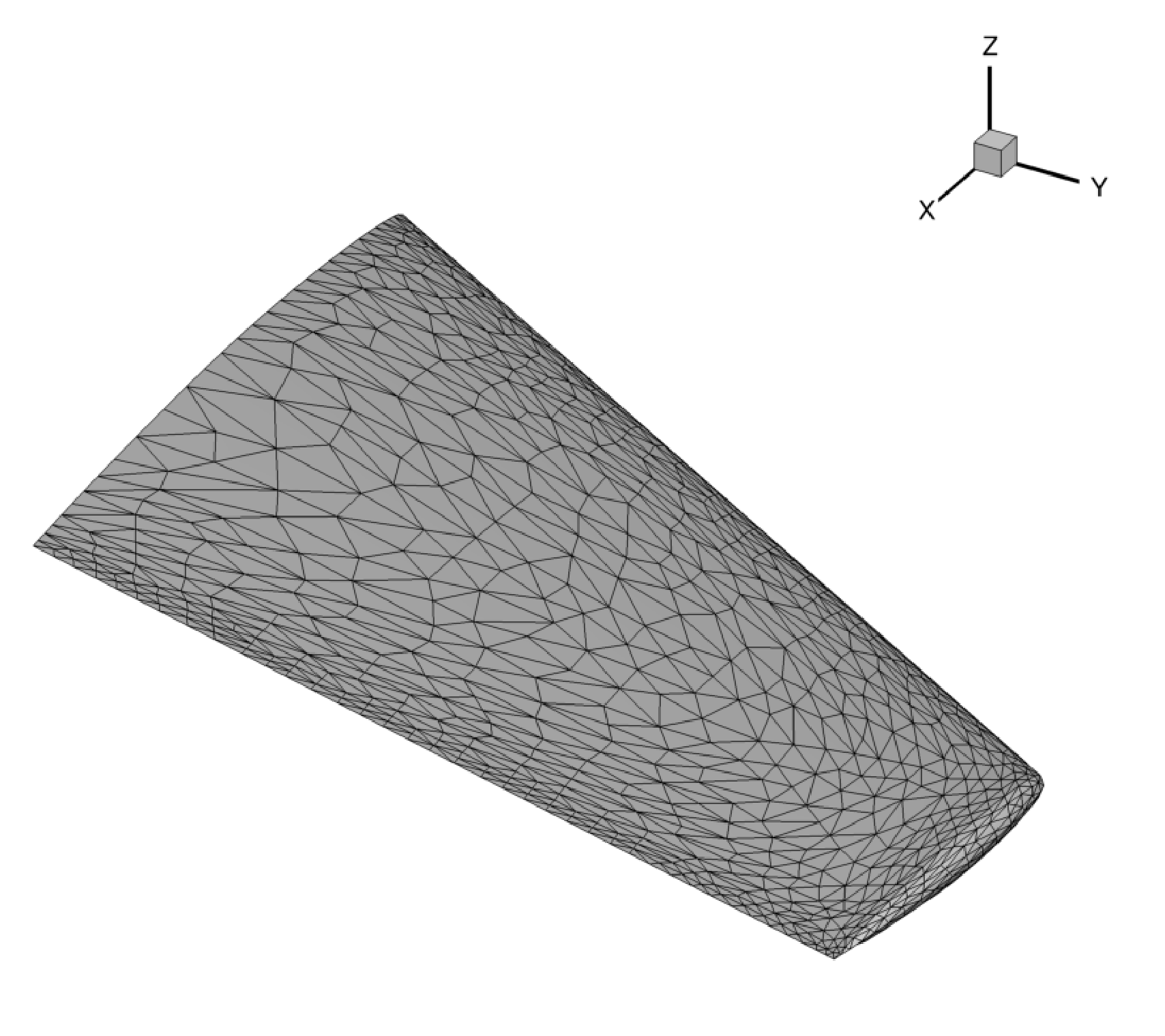}
    \caption{}\label{fig4:b}
    \end{subfigure}
    \begin{subfigure}[b]{.5\linewidth}
    \centering
    \includegraphics[width=.8\linewidth]{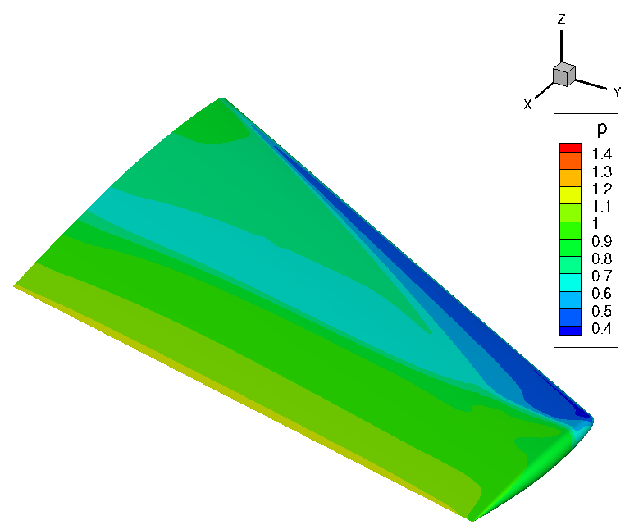}
    \caption{}\label{fig4:c}
    \end{subfigure}
    \begin{subfigure}[b]{.5\linewidth}
    \centering
    \includegraphics[width=.8\linewidth]{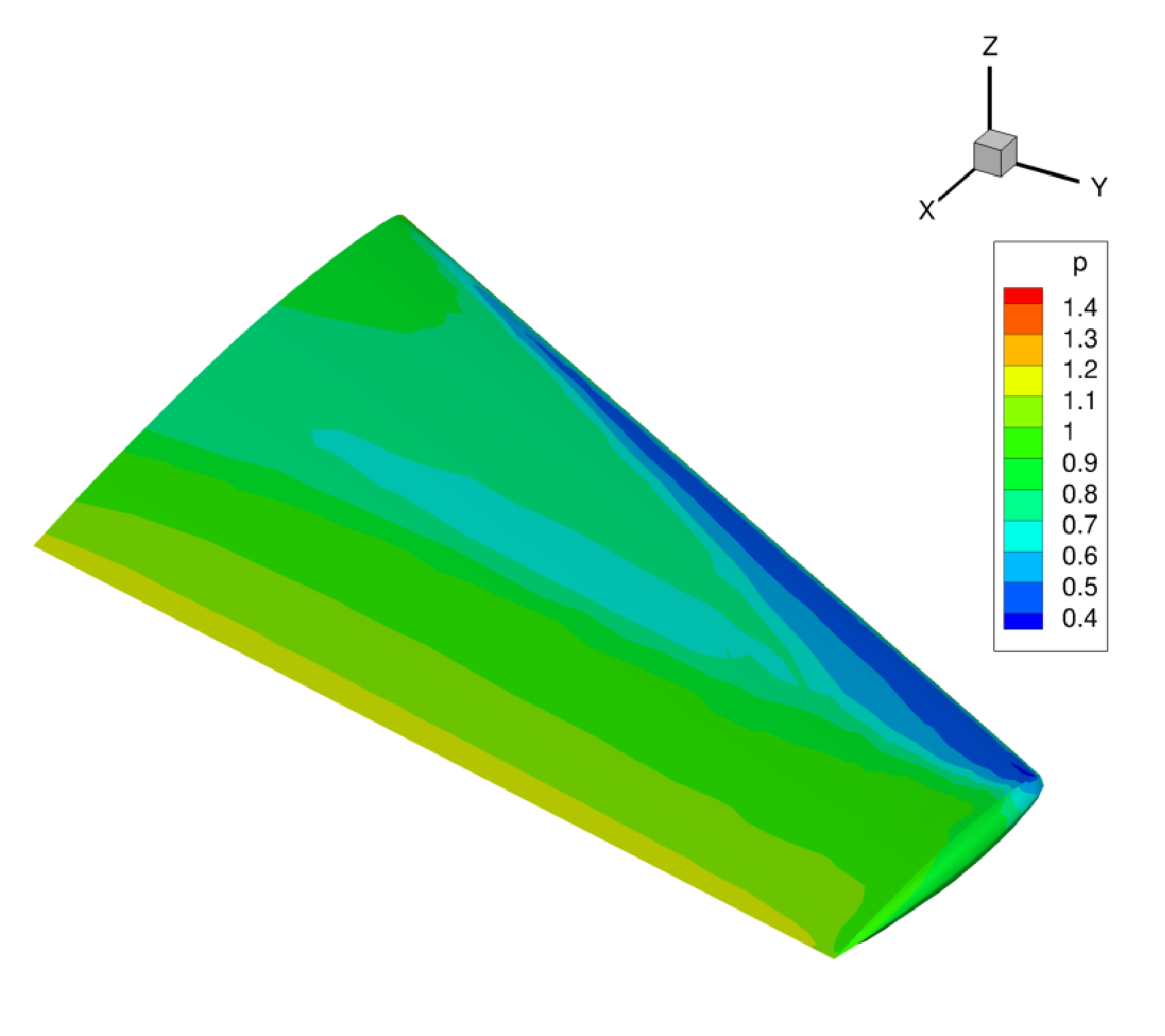}
    \caption{}\label{fig4:d}
    \end{subfigure}
    \caption{3D transonic turbulent flow over the ONERA M6 Wing: (a) and (b) represent respectively the wing geometry for structured and unstructured meshes while the steady state pressure field is illustrated in (c) and (d) for both FV and DG computations.}\label{fig4}
\end{figure}

The second test case corresponds to the steady transonic turbulent flow around the ONERA M6 wing. For the FV method, the structured mesh is composed of 60 blocks for a total of 3.6 million of cells. For the DG case, the unstructured high-order mesh is composed of curved quadratic tetrahedra and wedges for a total of 168 956 cells. Figure~\ref{fig4} depicts both the mesh and the steady state pressure field. For the FV scheme, the algebraic problem size is fixed to 21 million with a number of entries around 2 billion against 10 million with a number of entries around 3.2 billion for the DG scheme.

\subsubsection{Standard Krylov subspace methods}

This section aims at demonstrating how the standard GMRES solver suffers from numerical limitations when tackling stiff problems. It is crucial to point out this aspect since many authors frequently employ the standard GMRES in their numerical experiments. Hence, it is worthwhile to assess the behavior of standard GMRES compared to its flexible counterpart for stiff problems. To perform a fair comparison, the size of the respective Krylov subspaces are chosen in order to get the same memory footprint. Based on a maximum size of 60 and 20 vectors for the outer and inner Krylov spaces respectively, the memory footprint of FGMRES-DR amounts to the storage of 140 vector fields.

The relative true residual norm is usually expected to be of the same order to that for the convergence of the stationary problem.
But in several specific contexts, higher accurate solutions are required so we are looking for true residual norms as close as possible to machine epsilon.

\begin{figure}[ht]
    \centering
    \begin{subfigure}[b]{.5\linewidth}
    \centering
    \includegraphics[trim=0.25cm 0.5cm 1.5cm 1.5cm,clip,width=1.05\linewidth]{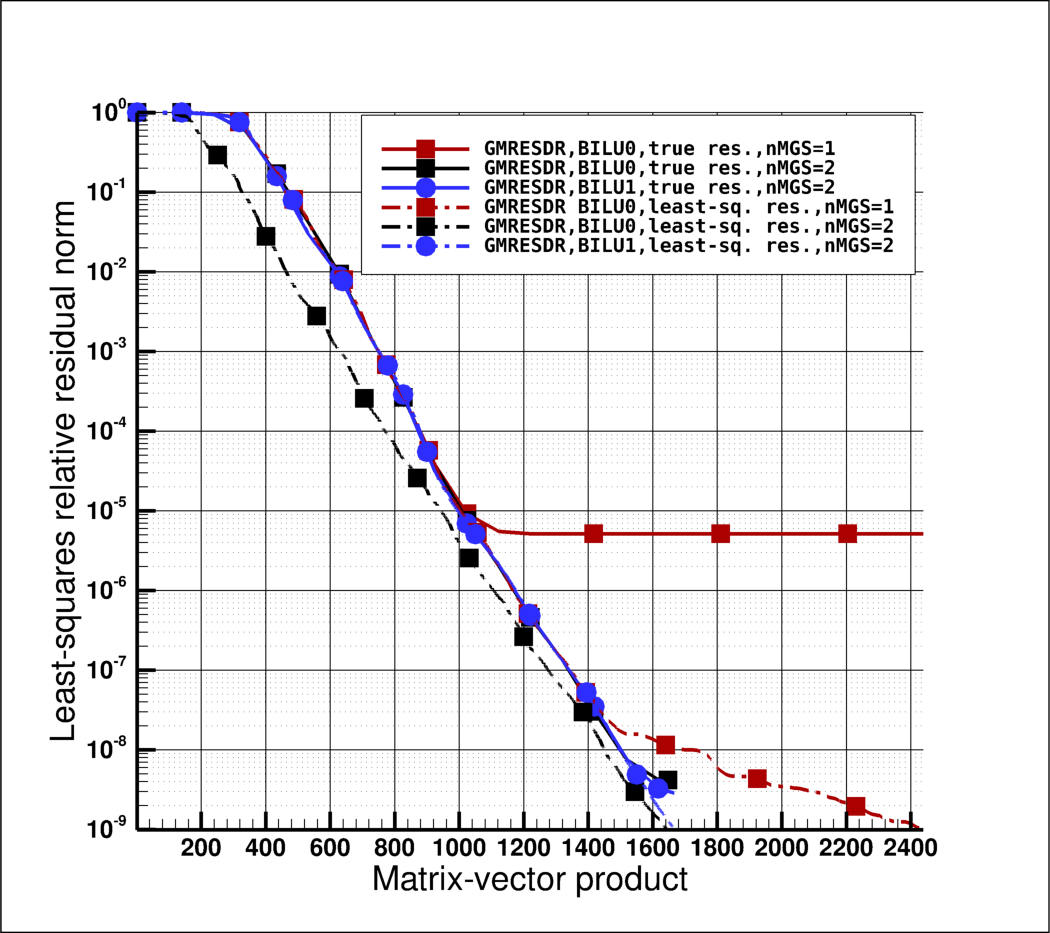}
    \caption{FV case}\label{fig5:a}
    \end{subfigure}
    \begin{subfigure}[b]{.5\linewidth}
    \centering
    \includegraphics[trim=0.25cm 0.5cm 1.5cm 1.5cm,clip,width=1.05\linewidth]{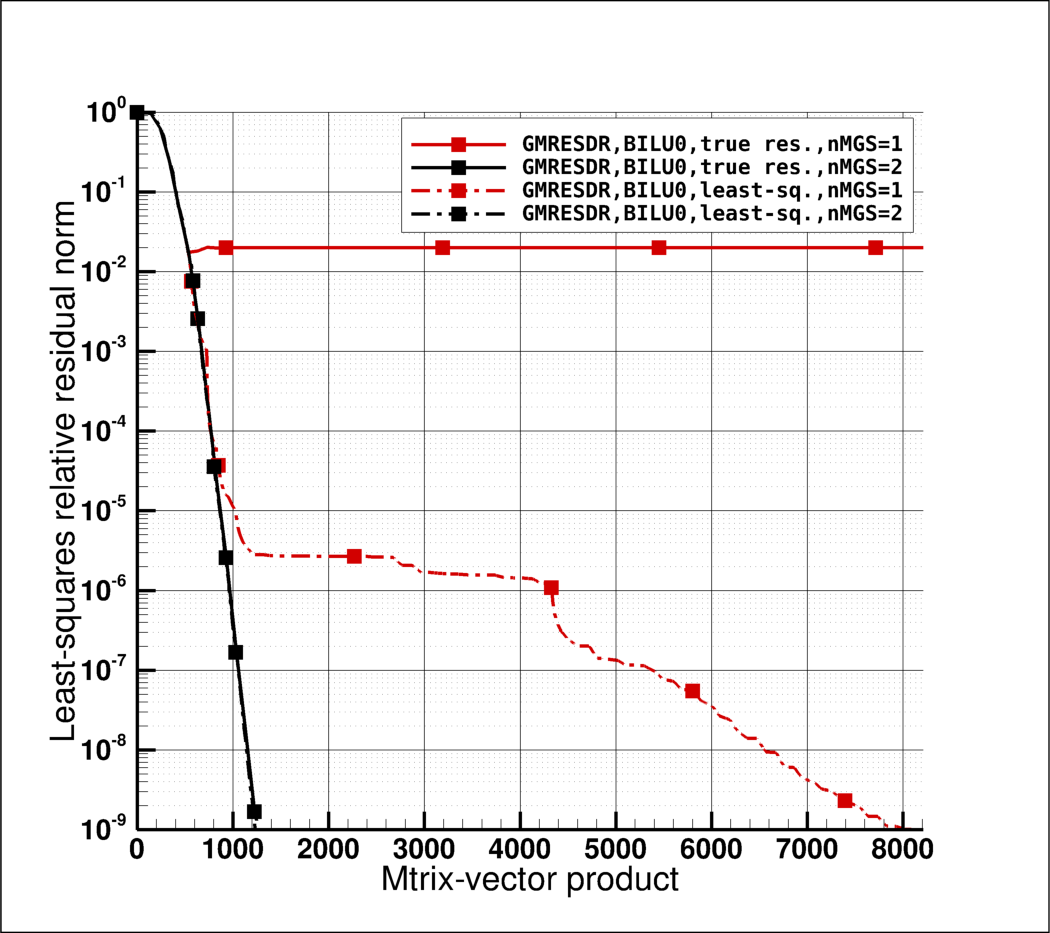}
    \caption{DG case}\label{fig5:b}
    \end{subfigure}
    \caption{True and least-squares residual norm convergence histories by applying two steps of Gram-Schmidt orthogonalization procedure: \emph{least-sq. res.} stands for least-squares residual while \emph{true res.} stands for true residual. We recall the numerical parameters of the GMRES-DR solver as follows: $m=140$ and $k=42$}\label{fig5}
\end{figure}

For such requirements, the Modified Gram-Schmidt process (MGS) exhibits limitations because it is known to perform poorly when input vectors are almost colinear. In this case a standard remedy is to apply a re-orthogonalization step thus doubling the cost of the Arnoldi procedure. For computational cost reasons, the residual convergence history is assessed with the least-squares residual while the quality of the convergence is given by the true residual which is only computed at the end of each cycle. Both residuals are plotted in Figure~\ref{fig5}.
All operations with the GMRES-DR solver are performed in a double-precision arithmetic.
When a single step of MGS orthogonalization is applied, the least-squares residual reaches the required tolerance while the corresponding true residual stagnates at $10^{-6}$ for BILU(0) applied to $\textbf{J}^{EX}_{O1}$ in the FV method. For the DG case, true and least-squares residuals are quite similar up to $10^{-2}$. At this level, the true residual stagnates while the least-squares residual decreases up to the prescribed tolerance.
This phenomenon comes essentially from the rounding errors during the orthogonalization process since the size of the Krylov subspace is sufficiently large to give rise to numerical errors. The corresponding convergence histories for two steps of orthogonalization are given in Figure~\ref{fig5}. For the FV problem, we get a gain of two orders of magnitude with a final relative true residual at $10^{-8}$. For the DG problem, the gain is even better with a converged relative true residual at $10^{-9}$. Another important aspect that explains why the GMRES-DR solver does not converge at the lower tolerance in the FV case comes from the quality of the preconditioner. We point out that the DG method uses essentially the third-order operator $\textbf{J}^{EX}_{O3}$ as a preconditioner, which corresponds to the same order of accuracy as the discretization scheme, while the FV method simply relies on approximate first order preconditioning operators while a second order spatial discretization scheme is applied for the mean flow solution. In an attempt to increase the quality of the preconditioner we used a BILU(1) factorization of on the first-order exact Jacobian matrix $\textbf{J}^{EX}_{O1}$. As can be seen in Figure~\ref{fig5}a, the corresponding gain is however marginal. 

\subsubsection{Inner-Outer Krylov subspace methods}

\begin{figure}[H]
    \centering
    \begin{subfigure}[b]{.49\linewidth}
    \centering
    \includegraphics[trim=0.25cm 0.5cm 1.5cm 1.5cm,clip,width=1.05\linewidth]{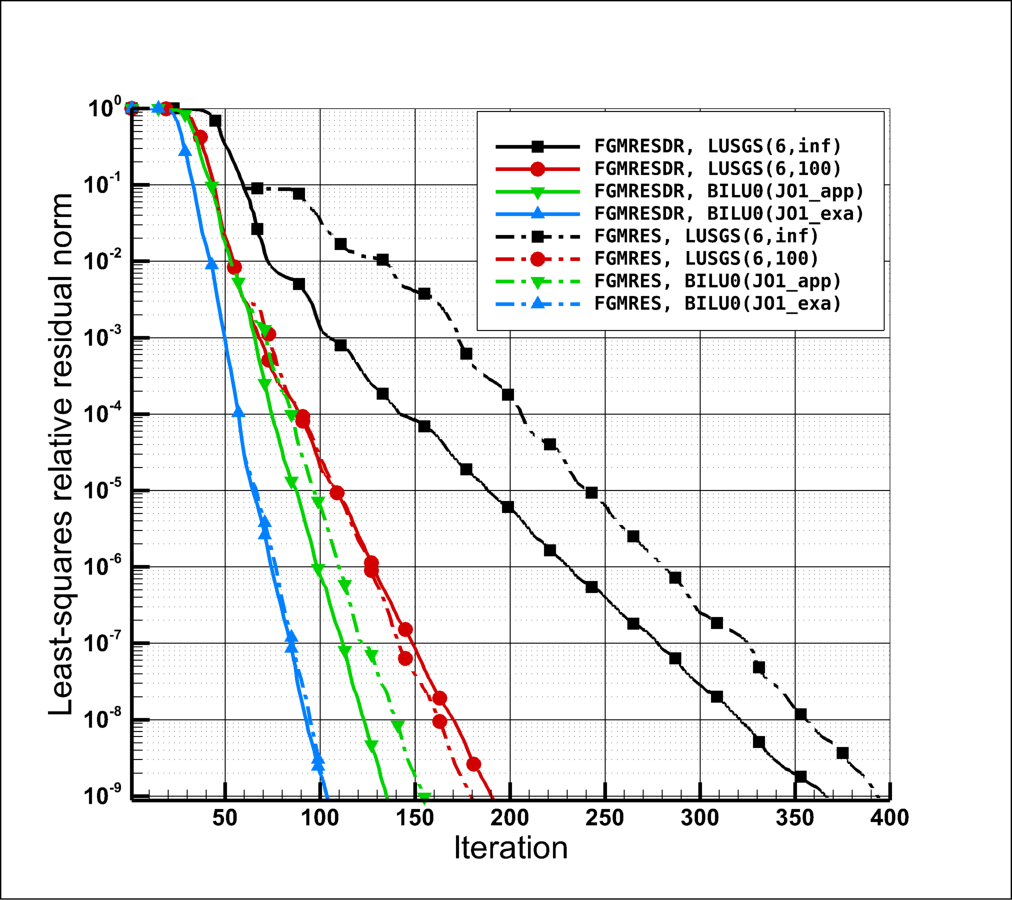}
    \caption{FV case}\label{fig6:a}
    \end{subfigure}
    \hfill
    \begin{subfigure}[b]{.49\linewidth}
    \centering
    \includegraphics[trim=0.25cm 0.5cm 1.5cm 1.5cm,clip,width=1.05\linewidth]{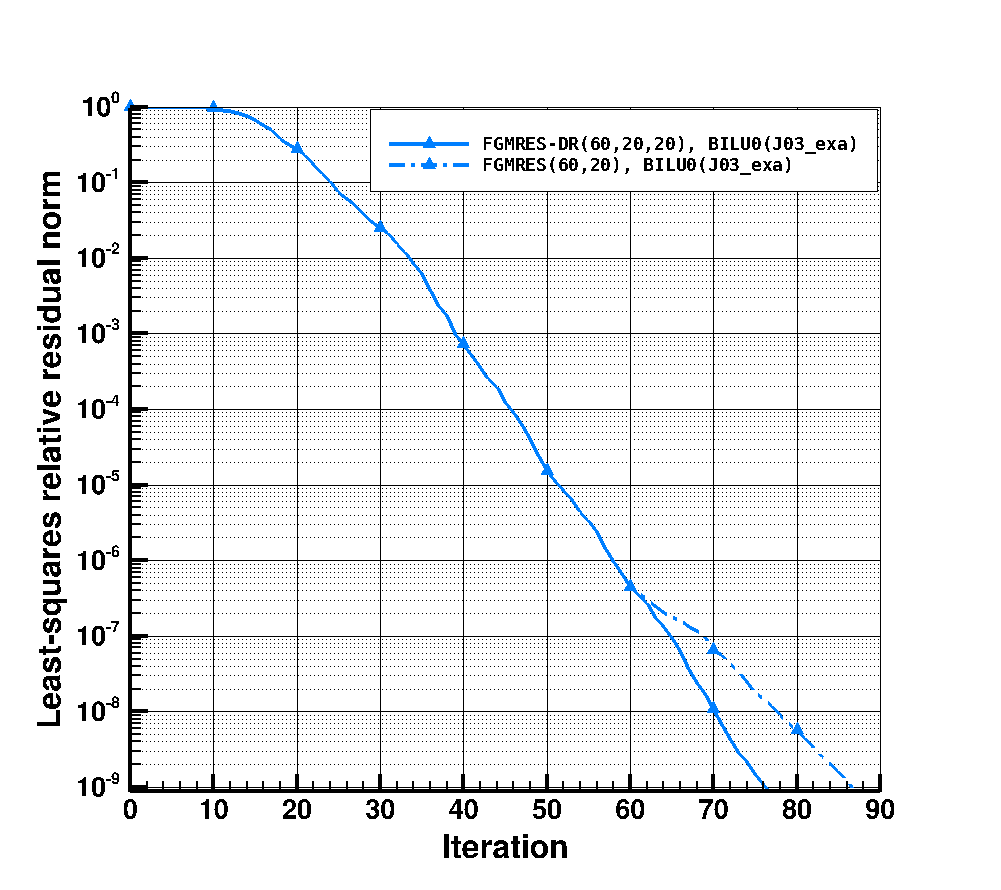}
    \caption{DG case}\label{fig6:b}
    \end{subfigure}
    \caption{Impact of deflation for various preconditioners. The relative residual norm convergence history is plotted with respect to iterations. The increase of the diagonal dominance of the LU-SGS preconditioner by addition of a scalar diagonal $I/{\Delta \tau}$ seems very promising but the strategy remains less effective than the one based on BILU. The impact of deflation is not very conclusive for the most efficient BILU preconditioner as only one restart is performed in this case. We recall the numerical parameters of the FGMRES-DR solver as follows: $m=60, m_{i}=20$ and $k=20$.}\label{fig6}
\end{figure}

For the nested strategy, we perform the same comparison study keeping the same numerical parameters detailed in Table~\ref{table1} of subsection~\ref{sec4_2}.
Unlike the two-dimensional OAT15A test case, the FGMRES(60,20) solver preconditioned either by LU-SGS(6) or BILU(0) applied to $\textbf{J}^{APP}_{O1}$ converges without the deflation strategy. In Figure~\ref{fig6}a, FGMRES(60,20) exhibits a good convergence with 150 iterations for BILU(0) against 377 for LU-SGS(6). However, the convergence is a little sensitive to the deflation effect with a gain of 22 iterations for LU-SGS(6) and 15 iterations for BILU(0). Even though the preconditioning matrix $\textbf{J}^{APP}_{O1}$ is less robust than $\textbf{J}^{EX}_{O1}$, it still remains interesting from a memory point of view as we have explained for the 2D test case. Actually, it is more crucial to take into account the memory factor when it comes to solve three-dimensional stiff problems. Indeed,  we recall that $\textbf{J}^{APP}_{O1}$ has a stencil of 7-points in 3D while a stencil of 19-points is associated with $\textbf{J}^{EX}_{O1}$, that is, a ratio of 2.7 against 1.8 for the two-dimensional case. 


\begin{figure}[H]
    \centering
    \begin{subfigure}[b]{.49\linewidth}
    \centering
    \includegraphics[trim=0.25cm 0.5cm 1.5cm 1.5cm,clip,width=1.05\linewidth]{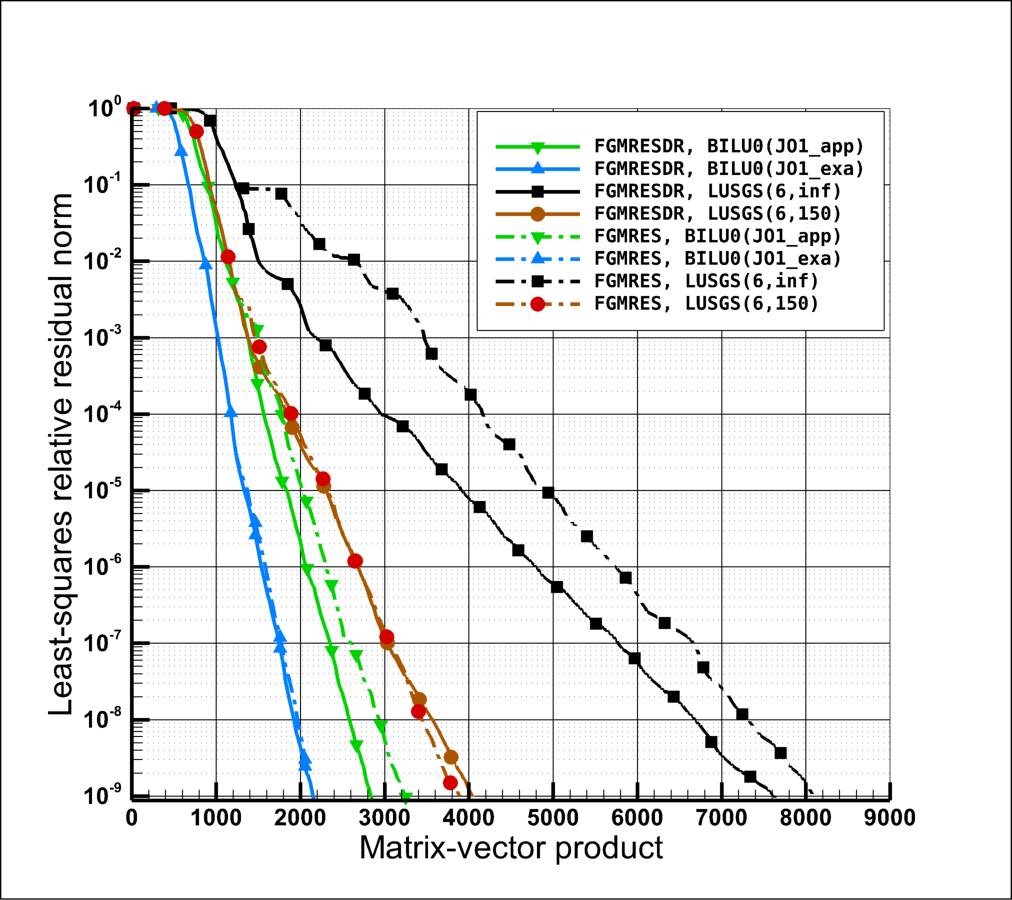}
    \caption{FV case}\label{fig7:a}
    \end{subfigure}
    \hfill
    \begin{subfigure}[b]{.49\linewidth}
    \centering
    \includegraphics[trim=0.25cm 0.5cm 1.5cm 1.5cm,clip,width=1.05\linewidth]{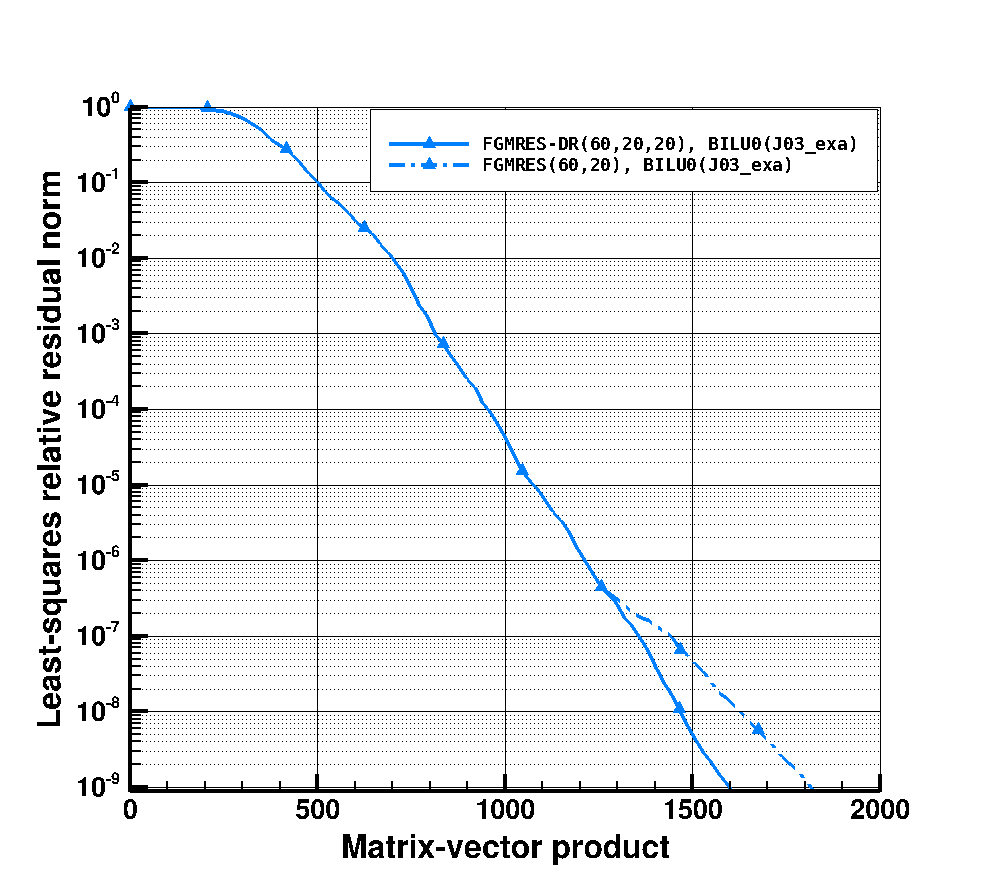}
    \caption{DG case}\label{fig7:b}
    \end{subfigure}
\caption{Impact of deflation for various preconditioners. The relative residual norm convergence history is plotted with respect to Jacobian-vector products. The hierarchy of the preconditioners remains the same as the one in Figure~\ref{fig6}.}\label{fig7}
\end{figure}

\begin{figure}[H]
   \begin{minipage}{0.49\textwidth}
     \centering
     \includegraphics[trim=0.25cm 0.5cm 1.5cm 1.5cm,clip,width=1.05\linewidth]{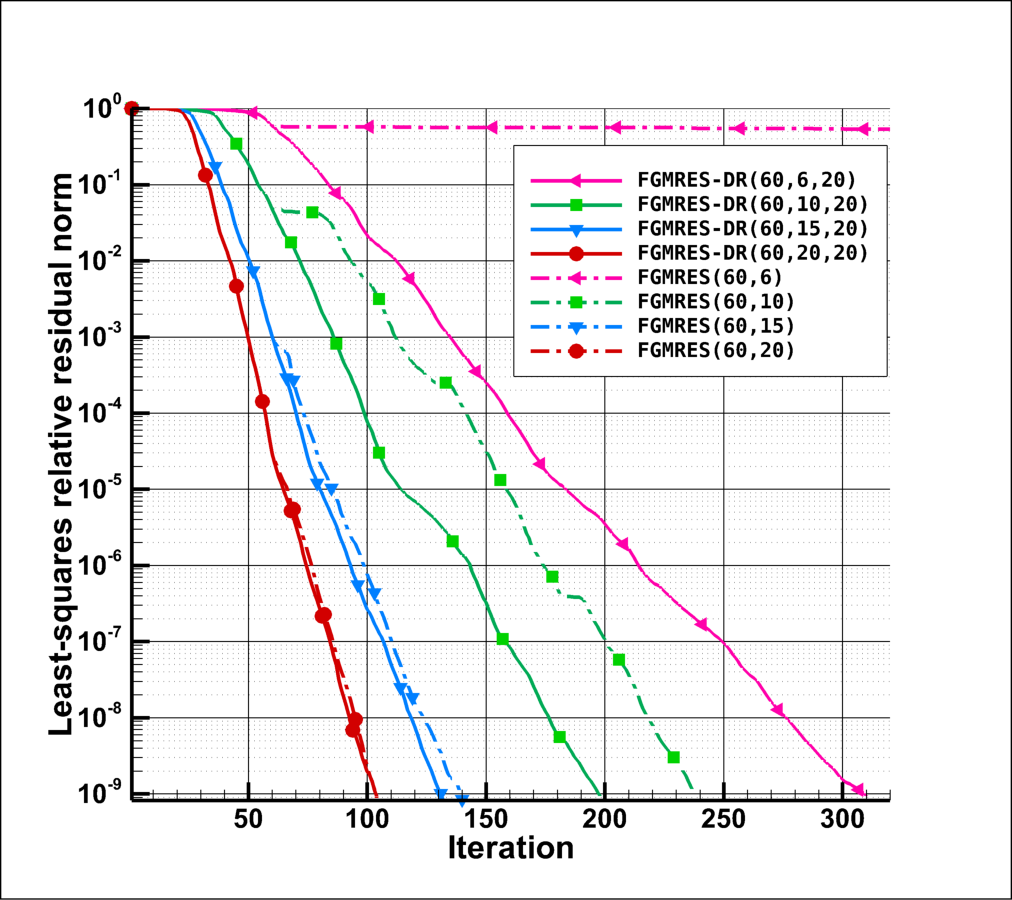}
   \end{minipage}\hfill
   \begin{minipage}{0.49\textwidth}
     \centering
     \includegraphics[trim=0.25cm 0.5cm 1.5cm 1.5cm,clip,width=1.05\linewidth]{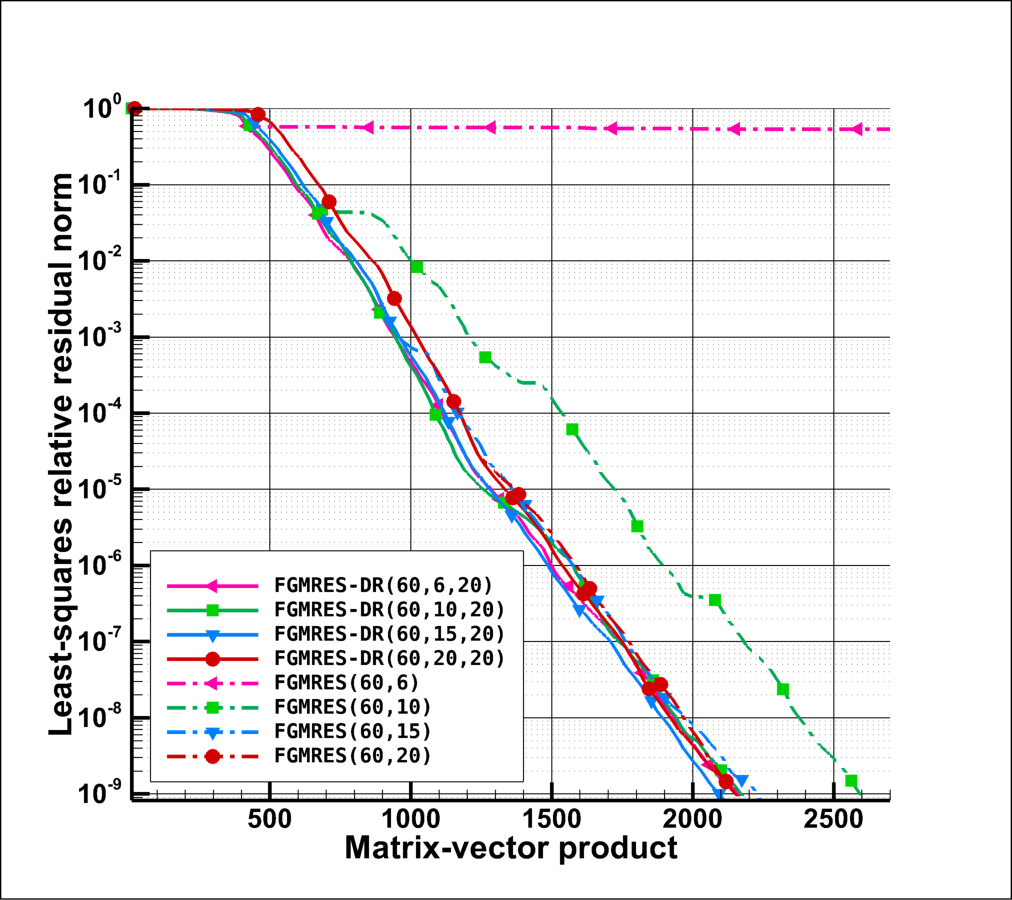}
   \end{minipage}
   \caption{Impact of deflation for various sizes of the inner Krylov space for FV case. Adjoint relative residual norm convergence history with respect to iterations (left) and Jacobian-vector products (right) for FGMRES($60,mi$) and FGMRES-DR($60,mi,20$) with $mi=6,10,15,20$.}\label{fig8}   
\end{figure}

In order to better highlight the impact of deflation we have conducted additional numerical simulations in the FV case by varying the size of the inner Krylov space with 6, 10, 15 and 20 vectors for the BILU type preconditioners. Figure~\ref{fig8} shows that when at least two restarts are performed the deflation step is clearly of interest in terms of iterations. Interestingly, the benefit in terms of matrix-vector products is not similar except for the lowest size of the inner Krylov space. This suggests that for this specific computation the best compromise in terms of memory and efficiency is to choose an inner Krylov size of 15.

\begin{figure}[H]
   \begin{minipage}{0.49\textwidth}
     \centering
     \includegraphics[trim=0.25cm 0.5cm 1.5cm 1.5cm,clip,width=1.05\linewidth]{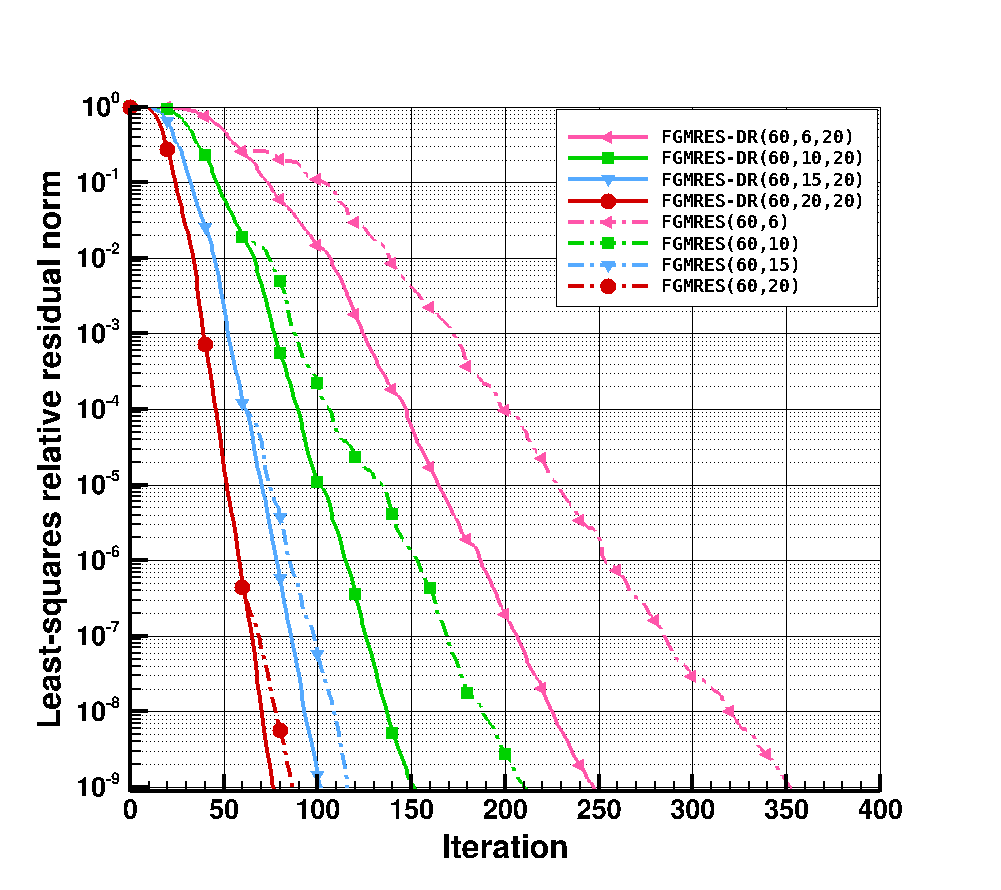}
   \end{minipage}\hfill
   \begin{minipage}{0.49\textwidth}
     \centering
     \includegraphics[trim=0.25cm 0.5cm 1.5cm 1.5cm,clip,width=1.05\linewidth]{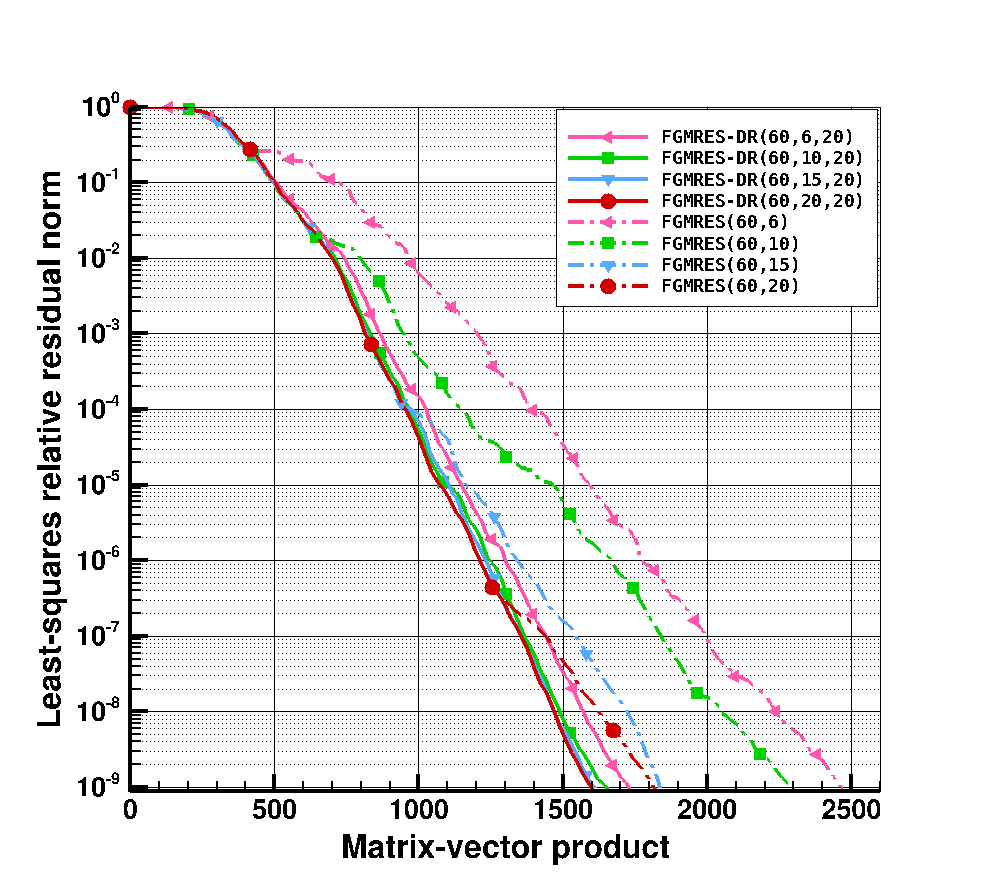}
   \end{minipage}
   \caption{Impact of deflation for various sizes of the inner Krylov space for DG case. Adjoint relative residual norm convergence history with respect to iterations (left) and Jacobian-vector products (right) for FGMRES($60,mi$) and FGMRES-DR($60,mi,20$) with $mi=6,10,15,20$.}\label{fig9}   
\end{figure}

\subsubsection{Scalability}

We have shown how the inner-outer GMRES solvers exhibit a satisfactory scalability with regard to the two-dimensional OAT15A airfoil. It would be advisable to retrieve the same scalability on the M6 wing configuration for both FV and DG formulations. 

For the scalability analysis, we consider five configurations corresponding to the partition of the initial domain into 2, 4, 8, 16 and 24 sub-domains to be performed on 24-core compute nodes.



For the unstructured case, the line construction algorithm~\cite{mavri2020} previously used for the two-dimensional OAT15A test-case has been applied. Finding a relevant set of weights for the METIS algorithm was more difficult on this three-dimensional test-case. A good trade-off has finally been found with a largest load imbalance about 5 $\%$. The number of iterations to converge the adjoint problem is sensitive to the number of cores with a variation up to 15 $\%$ compared to the two-subdomain reference configuration. As it could have a significant impact on such a challenging test-case, another strategy has been adopted and relies on the CHACO software~\cite{chaco1995} for partitioning graphs. Spectral methods partition the graph using the eigenvectors of a matrix constructed from the graph.
We select the simplest spectral method which is a weighted version of spectral bisection.
This method uses the second lowest eigenvector of the Laplacian matrix of the graph to divide it into two pieces~\cite{hdsimon1991}.
This results invariably to a quite perfect partitioning in terms of number of cells with almost one cell per partition of difference. The number of iterations to converge the adjoint problem shows now a variation up to 6 $\%$.

\begin{figure}[H]
    \centering
    \begin{subfigure}[b]{.5\linewidth}
    \centering
    \includegraphics[trim=0.2cm 0.5cm 1.5cm 1.5cm,clip,width=1.05\linewidth]{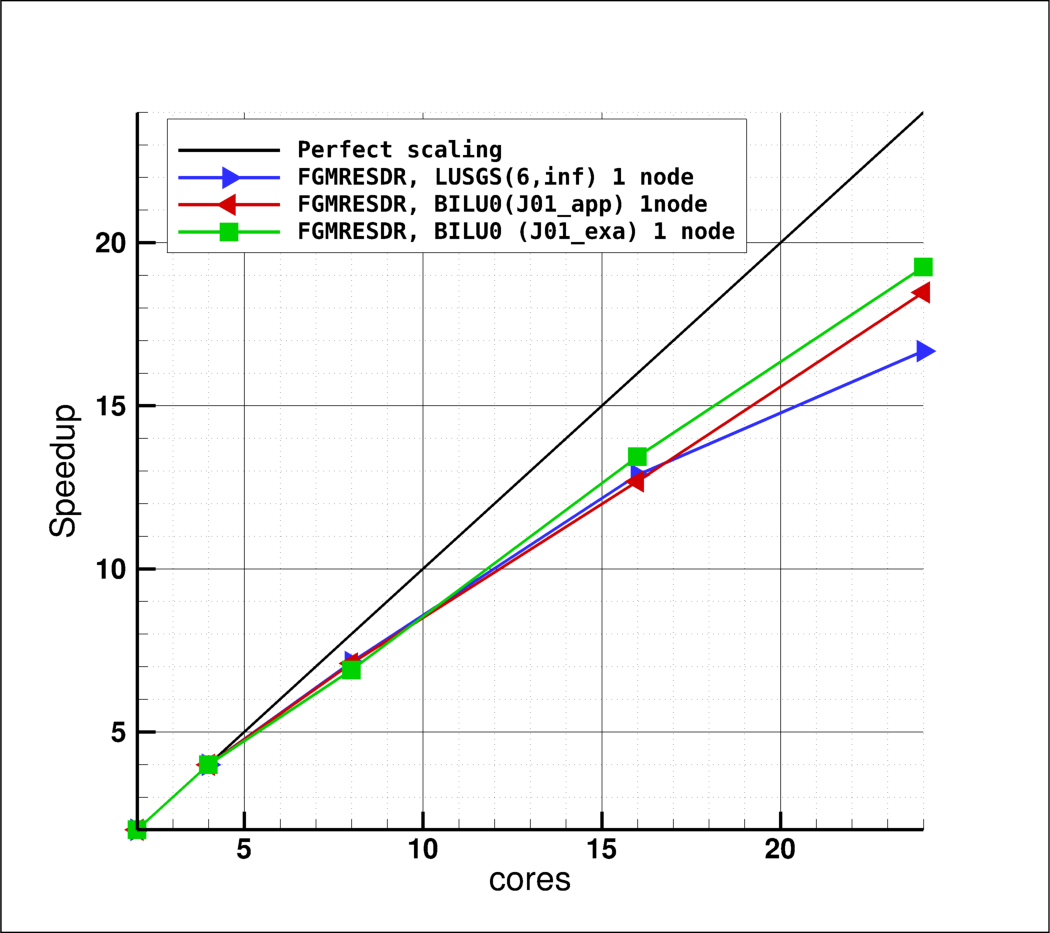}
    \caption{FV case}\label{fig10:a}
    \end{subfigure}
    \begin{subfigure}[b]{.5\linewidth}
    \centering
    \includegraphics[trim=0.2cm 0.5cm 1.5cm 1.5cm,clip,width=1.05\linewidth]{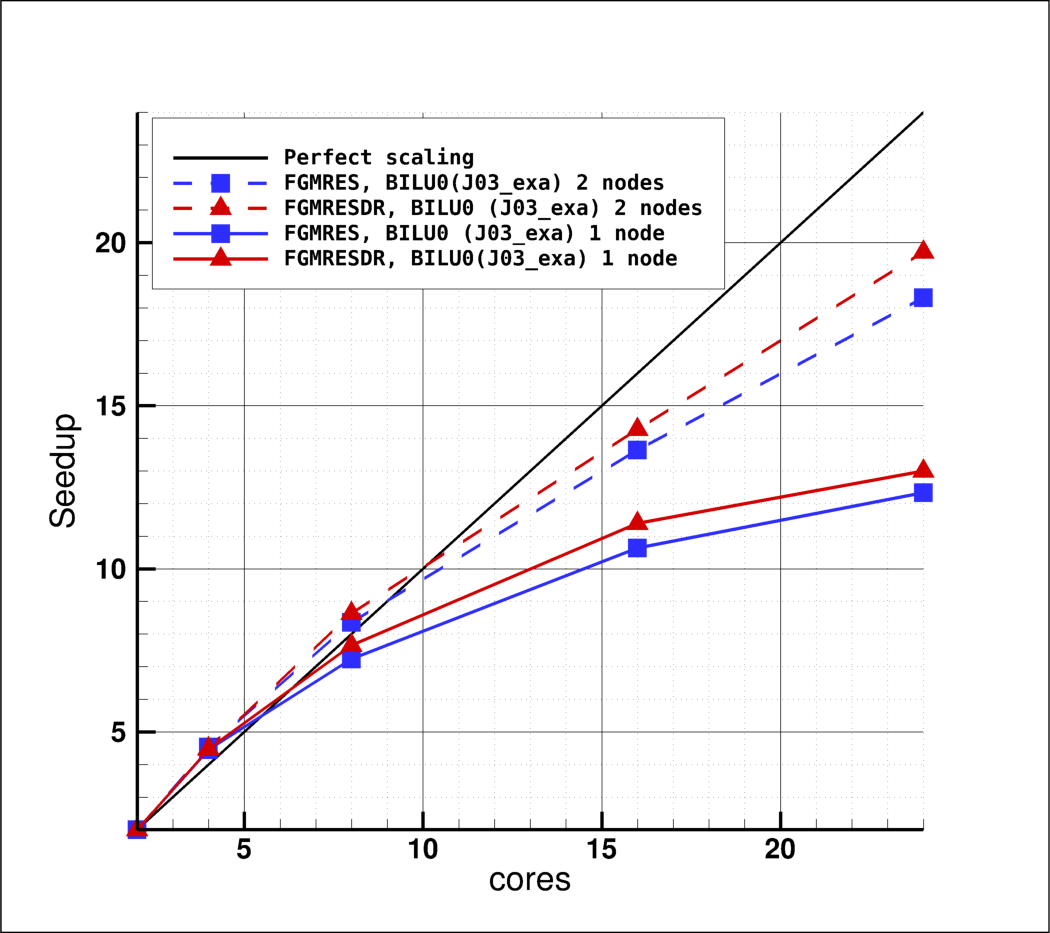}
    \caption{DG case}\label{fig10:b}
    \end{subfigure}
    \caption{Strong scalability analysis for FGMRES-DR(60,20,20) solvers on M6 wing adjoint system: both 1 node and 2 nodes have been plotted in the DG case.}\label{fig10}
\end{figure}

Figure~\ref{fig10} depicts parallel behavior of the FGMRES-DR solver. In Figure~\ref{fig10:b}, for the DG case, scalability analysis is performed on one node with one MPI process per core (solid line), and on two nodes (dashed line), with a balanced distribution of the MPI processes per node. This exercise exhibits a memory-bandwidth limitation for this test-case : the parallel efficiency is about 55$\%$ on one node and grows up to 85 $\%$ on two nodes. Dense linear algebra kernels (BLAS, LAPACK) are applied on blocks of size $60 \times 60$, corresponding to the number of degrees of freedom per mesh element, of the implicit matrix. But there are memory-bound operations with a low arithmetic intensity.

For the FV case, the same exercise has been performed with different scalability results. More specifically, the parallel efficiency is about 80 $\%$ regarding the BILU(0) preconditioner and 70 $\%$ for the LU-SGS one on one node. When it comes to use two nodes, the parallel efficiency grows up to 80 $\%$ for all preconditioners. The matrix operations now work on blocks of smaller size $6\times6$ and are therefore less prone to memory-bandwith limitations.\\

Since the preconditioner does depend on the number of sub-domains, we have reported in Table~\ref{table4} the variation of the number of iterations of the FGMRES-DR according to the number of the sub-domains.

\begin{table}[H]
\centering
\begin{tabular}{|l|*{5}{c|}}\hline
\backslashbox{ $\#$ its ($\mathcal{M}$)}{$\#$ sub-domains}
&\makebox[3em]{2}&\makebox[3em]{4}&\makebox[3em]{8}
&\makebox[3em]{16}&\makebox[3em]{24}\\ \hline
$\#$ its (LUSGS(6,150)) & 190 & 189 & 190 & 190 & 190\\\hline
$\#$ its (LUSGS(6,inf)) & 348 & 348 & 349 & 349 & 348\\\hline
$\#$ its (BILU(0), $\textbf{J}^{O1}_{APP}$) & 134 & 135 & 134 & 134 & 135\\\hline
$\#$ its (BILU(0), $\textbf{J}^{O1}_{EX}$) & 99 & 98 & 99 & 99 & 98\\\hline
\end{tabular}
\caption{Impact of the number of sub-domains on the number of iterations for FGMRES-DR(60,20,20). As the inner GMRES preconditionner is global by construction, we see no impact on the number of iterations as expected.}
\label{table4}
\end{table}

\section{Some numerical considerations}

\hspace{0.5cm} In our numerical experiments we have considered low convergence thresholds because it is acknowledged that adjoint solutions may exhibit oscillations or still be locally under-resolved for moderate relative residual decrease (typically 3 orders of magnitude). This strong convergence is only attainable in practice if a number of precautions is taken in the implementation of the FGMRES-DR solver. First, the orthogonality of the Krylov basis at the outcome of the Arnoldi process or after a deflation step must be preserved as much as possible to ensure convergence in a reasonable number of iterations for stiff problems. This point was studied in~\cite{rollin2008improving} where the authors advocate the use of a re-orthogonalization during the MGS process. Morgan~\cite{morgan2002} also proposed an economical strategy where only the last Krylov vector of the deflated basis is re-orthogonalized with respect to other deflated basis vectors. The initial physical equations should always be solved in their non dimensionalized form. For large ill-conditioned linear systems, Saad~\cite{saad2003} also recommends to apply a scaling to all the rows (or columns) e.g., so that their 1-norms are all equal to 1; then to apply a scaling of the columns (or rows). This practice is strongly recommended before any ILU type factorization. Additionally, recognizing that the linearization of the turbulence model dramatically deteriorates the condition number of the system matrix, Chisholm and Zingg~\cite{chisholm2009} advice to apply a specific scaling coefficient to this particular equation. Additional numerical ingredients like the mixed-precision nested GMRES, i.e., single precision floating point arithmetic for the inner GMRES, and the application of the BILU factorization to the RAS partition of the fluid domain have proven beneficial. Last but not least, the availability of different physical approximations to build the candidate flux Jacobian matrix for the preconditioner is very interesting to keep a good compromise between efficiency and memory footprint as it has a strong impact on the size of the stencil of this operator. Finally, some memory saving improvements were proposed in~\cite{giraud2010} to perform inplace the matrix-matrix product leading to the new deflated Krylov vector basis at the starting of the next cycle. In fact in a nested FGMRES-DR framework this is not required anymore because the temporary memory required by this product is already available for the storage of the inner Krylov subspace and is then naturally re-used during the outer level deflation step. It is also worth mentioning that in practice the inner GMRES solver is never restarted considering the small size of the inner Krylov subspace. 

Another important aspect that should be highlighted is the numerical limitation of the Modified Gram-Schmidt process. 
Indeed, in terms of number of domains in the mesh partition, it requires asymptotically more messages and synchronization points in parallel, and asymptotically more data traffic between levels of the memory hierarchy.

\section{Concluding remarks}

Flexible GMRES solvers have been investigated for the solution of large
sparse linear systems arising from discrete stiff adjoint problems when
the standard GMRES fails. A FV and a high-order DG discretizations are
considered to provide linear systems of different complexity from
test-cases of optimization problems of turbulent flows with RANS
modelling. Condition number and dimensions of the resulting matrices
change but also performance of the underlying basic linear algebra
operations. Varying the discretization methods contributes to a better
overview of the solver capabilities. Relevant numerical ingredients have
been developed to improve robustness and efficiency of the inner-outer
GMRES strategy and concern preconditioning techniques and deflation
strategy. Numerical experiments have been conducted on turbulent
transonic flows over the two-dimensional supercritical ONERA OAT15A
airfoil and over the three-dimensional ONERA M6 wing.

For the first ingredient, with the FV scheme, we have investigated the
impact of using different approximations of the flux Jacobian matrix for
the preconditioning step. The efficiency of the nested GMRES has been
observed when deflation technique with the first order approximate
Jacobian matrix is considered or by using the first order exact Jacobian
matrix without deflation but at the price of twice the memory cost. The
use of an exact linearization of the DG scheme to build the
preconditioner offers very good performance of the linear solvers in
terms of robustness and convergence.

For the second ingredient, significant gains has been obtained on the
convergence speed. Not only eigen-information can be recovered at low
cost using by-products of the Arnoldi iterative process, but
matrix-vector products involving the Jacobian matrix can be saved at the
beginning of each new cycle thanks to the recycling of eigenvectors into
the new Krylov basis.

On top of that, a good parallel scalability of the FGMRES is achieved
thanks to the inherent global effect of the inner GMRES. The interesting
conclusion of this work is the capability of inner-outer GMRES strategy
to solve stiff problems and exhibit the superlinear convergence of the
full GMRES, but with a lower memory footprint. Strategies to select
eigen-information and calibration rules for the numerical parameters
driving the solution of inner and outer systems might be worth taking
into consideration in the future.

\newpage

\bibliography{CAF-D-21-00578}

\end{document}